\documentclass[a4paper,10pt]{article}
\usepackage{graphicx}
\usepackage{amsmath}
\usepackage{amsfonts}
\usepackage{amssymb}
\usepackage{amsthm}

\newtheorem{theorem}{Theorem}[section]
\newtheorem{lemma}{Lemma}[section]
\newtheorem{remark}{Remark}[section]
\newtheorem{proposition}{Proposition}[section]

\newtheorem{corollary}{Corollary}[section]
\newtheorem{algorithm}{Algorithm}[section]
\renewcommand{\theequation}{\arabic{section}.\arabic{equation}}
\renewcommand{\thetheorem}{\arabic{section}.\arabic{theorem}}
\renewcommand{\thelemma}{\arabic{section}.\arabic{lemma}}
\renewcommand{\theproposition}{\arabic{section}.\arabic{proposition}}

\begin{document}
\title{Convergence of Adaptive Finite Element Approximations
for Nonlinear Eigenvalue Problems \thanks{{This work was partially
supported by the National Science Foundation of China under grants
10871198 and 10971059, the National Basic Research Program of China
under grant 2005CB321704, and the National High-Tech Development
Program of China under grant 2009AA01A134.}}}

\author{Huajie Chen\thanks{LSEC, Institute of Computational Mathematics
and Scientific/Engineering Computing, Academy of Mathematics and
Systems Science, Chinese Academy of Sciences, Beijing 100190, China
({hjchen@lsec.cc.ac.cn}).} \and Xingao Gong\thanks{ Department of
Physics, Fudan University, Shanghai 200433, China
(xggong@fudan.edu.cn).} \and  Lianhua He\thanks{LSEC, Institute of
Computational Mathematics and Scientific/Engineering Computing,
Academy of Mathematics and Systems Science, Chinese Academy of
Sciences, Beijing 100190, China ({helh@lsec.cc.ac.cn}).}\and Aihui
Zhou\thanks{LSEC, Institute of Computational Mathematics and
Scientific/Engineering Computing, Academy of Mathematics and Systems
Science, Chinese Academy of Sciences, Beijing 100190, China
({azhou@lsec.cc.ac.cn}).}}

\date{}
\maketitle

\begin{abstract}
In this paper, we study an adaptive finite element method for a
class of a nonlinear eigenvalue problems that may be of nonconvex
energy functional and consider its applications to quantum
chemistry. We prove the convergence of adaptive finite element
approximations and present several numerical examples of
micro-structure of matter calculations that support our theory.
\end{abstract}



\section{Introduction}\setcounter{equation}{0}
In this paper, we study adaptive finite element approximations for a
class of nonlinear eigenvalue problems: Find $\lambda\in \mathbb{R}$
and $u\in H^1_0(\Omega)$ such that
\begin{eqnarray}\label{problem-eigen}
\left\{ \begin{array}{rcl} \left(-\alpha\Delta + V +
\mathcal{N}(u^2) \right) u &=& \lambda u \quad \mbox{in} \quad
\Omega,
\\[1ex] \int_{\Omega}|u|^2 &=& Z,
\end{array} \right.
\end{eqnarray}
where $\Omega\subset \mathbb{R}^3$, $Z\in \mathbb{N},~\alpha\in
(0,\infty)$, $V:\Omega\rightarrow \mathbb{R}$ is a given function,
$\mathcal{N}$ maps a nonnegative function over $\Omega$ to some
function defined on $\Omega$.

Many physical models for micro-structures of matter are nonlinear
eigenvalue problems of type  (\ref{problem-eigen}), for instance,
the Thomas-Fermi-von Weizs\"{a}cker (TFW) type orbital-free model
used for electronic structure calculations
\cite{chen-zhou08,lieb81,wang-carter00} and the Gross-Pitaevskii
equation (GPE) describing the Bose-Einstein condensates (BEC)
\cite{bao-du04,zhou04}. In the context of simulations of electronic
structure calculations, the basis functions used to discretize
(\ref{problem-eigen}) are traditionally plane wave bases or
typically Gaussian approximations of the eigenfunctions of a
hydrogen-like operator. The former is very well adapted to solid
state calculations and the latter is incredibly efficient for
calculations of molecular systems. However, there are several
disadvantages and limitations involved in such methods. For example,
the boundary condition does not correspond to that of an actual
system;  extensive global communications in dealing with plane waves
reduce the efficiency of a massive parallelization, which is
necessary for complex systems; and the generation of a large
supercell is needed for non-periodic systems, which certainly
increases the computational cost. The finite element method uses
local piecewise polynomial basis functions, which does not involve
problems mentioned above and has several advantages. Although it
uses more degrees of freedom than that of  traditional methods,
strictly local basis functions produce well structured sparse
Hamiltonian matrices; arbitrary boundary conditions can be easily
incorporated; more importantly, since ground state solutions
oscillate obviously near the nuclei, it is relatively
straightforward to implement adaptive refinement techniques for
describing regions around nuclei or chemical bonds where the
electron density varies rapidly, while treating the other zones with
a coarser description, by which computational accuracy and
efficiency can be well controlled. Thus it should be  natural to
apply adaptive finite element methods to solve nonlinear eigenvalue
problems resulting from modeling electronic structures. Indeed the
adaptive finite element method is a powerful approach to computing
ground state energies and densities in quantum chemistry, materials
science, molecular biology and nanosciences \cite{beck00,bris03}.

The basic idea of a standard adaptive finite element method  is to
repeat the following procedure until a certain accuracy is obtained:
$$
\mbox{Solve}~\rightarrow~\mbox{Estimate}~\rightarrow~
\mbox{Mark}~\rightarrow~\mbox{Refine}.
$$
Adaptive finite element methods have been studied extensively since
Babu{\v s}ka and Rheinboldt \cite{babuska-rheinboldt78} and have
been successful in the practice of engineering and scientific
computing. In particular, D{\" o}rfler \cite{dorfler96} presented
the first multidimensional convergence result, which has been
improved and generalized, see, e.g.,
\cite{binev-dehmen-devore04,cascon-kreuzer-nochetto-siebert08,
mekchay-nochetto05,morin-nochetto-siebert00,
morin-nochetto-siebert02,morin-siebrt-veeser08, stevenson06a} for
linear boundary value problems,
\cite{carstensen08,chen-holst-xu07,dofer95,he-zhou09,veeser02} for
nonlinear boundary value problems, and
\cite{carstensen09,dai-xu-zhou08,garau-morin-zuppa09,garau-morin09,Giani-Graham09}
for linear eigenvalue problems. To our best knowledge,  there has
been no  work on the convergence of adaptive finite element
approximations for nonlinear eigenvalue problems, though some a
priori error analyses of finite dimensional Galerkin discretizations
for such nonlinear eigenvalue problems have been shown in
\cite{cances-chakir-maday09a, cances-chakir-maday09b,
chen-gong-zhou09,langwallner-ortner-suli09,zhou04,zhou07}.

In this paper, we shall present a posteriori error analysis of an
adaptive finite element method for a class of nonlinear eigenvalue
problems and prove that the adaptive finite element algorithm will
produce a sequence of approximations that converge to  exact ground
state solutions.  As an illustration, we shall also report several
numerical experiments on electronic structure calculations based on
the adaptive finite element discretization
\cite{beck00,bris03,chen-zhou08}, which support our theory. Since
the nonlinear term occurs, especially the nonlocal convolution
integration part, there are several serious difficulties in the
numerical analysis. Moreover, the associated energy functional for
this type of problems is usually nonconvex, which brings serious
difficulties.  In our analysis, we shall apply some nonlinear
functional arguments and special techniques to deal with the local
and nonlocal terms carefully.

 This paper is organized as
follows. In the coming section, we give an overview of the nonlinear
eigenvalue problem. In Section \ref{sec-FED}, we describe the finite
element discretization and give an a posteriori error analysis. In
Section \ref{sec-convergence},  we design an adaptive finite element
algorithm and prove  the convergence of  the adaptive finite element
algorithm. In Section \ref{sec-example}, we show some numerical
results for micro-structure computations to support our theory.
Finally, we give several concluding remarks.

\section{Preliminaries}\label{sec-preliminaries}
\setcounter{equation}{0} Let $\Omega\subset\mathbb{R}^3$ be a
polytypic bounded domain. We shall use the standard notation for
Sobolev spaces $W^{s,p}(\Omega)$ and their associated norms and
seminorms, see, e.g., \cite{adams75,ciarlet78}. For $p=2$, we denote
$H^s(\Omega)=W^{s,2}(\Omega)$ and $H^1_0(\Omega)=\{v\in H^1(\Omega):
v\mid_{\partial\Omega}=0\}$, where $v\mid_{\partial\Omega}=0$ is
understood in the sense of trace, $\|\cdot\|_{s,\Omega}=
\|\cdot\|_{s,2,\Omega}$, and $(\cdot,\cdot)$ is the standard $L^2$
inner product. The space $H^{-1}(\Omega)$, the dual of
$H^1_0(\Omega)$, will also be used. For convenience, the symbol
$\lesssim$ will be used in this paper. The notation $A\lesssim B$
means that $A\leq C B$ for some constant $C$ that is independent of
mesh parameters. We shall use $Pol(p,(c_1,c_2))$ to denote a class
of functions that satisfy the growth condition:
\begin{align*}\label{pol-notation}
Pol(p,(c_1,c_2))=&~\big\{ f:~ \exists ~a_1,a_2\in \mathbb{R} \mbox{
such that } \\& ~~c_1 |t|^{p}+a_1 \leq f(t) \leq c_2 |t|^{p}+a_2
\quad \forall ~t\ge 0 \big\}
\end{align*}
with $c_1\in \mathbb{R}$ and $c_2, p\in [0,\infty)$.

The weak form of (\ref{problem-eigen}) reads as follows: Find
$\lambda\in\mathbb{R}$ and $u\in H^1_0(\Omega)$ such that
\begin{eqnarray}\label{weak}
\left\{ \begin{array}{rcl} \alpha(\nabla u,\nabla v)
+(Vu+\mathcal{N}(u^2)u,v) & = & \lambda(u,v)
\quad \forall~ v\in H^1_0(\Omega), \\[1ex]
\|u\|^2_{0,\Omega} &=& Z.\end{array} \right.
\end{eqnarray}
For convenience, we divide  nonlinear term $\mathcal{N}$  into local
and nonlocal parts:
\begin{eqnarray*}
\mathcal{N}(\rho)=\mathcal{N}_1(\rho)+\mathcal{N}_2({\rho}),
\end{eqnarray*}
where $\rho=u^2$, $\mathcal{N}_1:[0,\infty)\rightarrow\mathbb{R}$ is
a given function dominated by some polynomial, and $\mathcal{N}_2$
is represented by a convolution integration
\begin{eqnarray*}\label{nonlinearepart2-}
\mathcal{N}_2(\rho) = \rho^{q-1}\int_{\Omega}\rho^{q}(y)K(\cdot-y)dy
\end{eqnarray*}
for some given function $K$ and $q\in\mathbb{R}$.

The associated energy functional with respect to this nonlinear
eigenvalue problem is expressed by
\begin{eqnarray}\label{energy}
E(u)= \int_{\Omega} \left(\alpha|\nabla u(x)|^2 + V(x)u^2(x) +
\mathcal{E}(u^2(x))\right)dx + \frac{1}{2q} D_K(u^{2q},u^{2q}),
\end{eqnarray}
where $\mathcal{E} :[0,\infty)\rightarrow \mathbb{R}$ is associated
with $\mathcal{N}_1$:
\begin{eqnarray*}\label{F-derive}
\mathcal{E}(s)=\int_{0}^{s}\mathcal{N}_1(t)dt,
\end{eqnarray*}
and $D_K(\cdot,\cdot)$ is a bilinear form defined by
\begin{eqnarray*}\label{D-convolution}
D_K(f,g) = \int_{\Omega} \int_{\Omega} f(x)g(y)K(x-y)dxdy.
\end{eqnarray*}
The ground state solution of problem (\ref{problem-eigen}) is
obtained by minimizing  energy functional (\ref{energy}) in the
admissible class
\begin{eqnarray*}\label{admissible-class}
\mathcal{A} = \left\{ \psi\in H^1_0(\Omega):
\quad\|\psi\|^2_{0,\Omega}=Z,\quad \psi\geq 0 \right\}.
\end{eqnarray*}

In our discussion, we assume that
\begin{itemize}
\item[(i)] $V\in L^2(\Omega)$.

\item[(ii)] $\mathcal{E}\in Pol(p,(c_1,c_2))$ with one of the
following conditions:
\begin{itemize}
\item[1.] $c_1\in(0,\infty)$;
\item[2.] $p\in [0,4/3]$;
\item[3.] $c_1\in
(-\infty,0)$, $p\in(4/3,\infty)$ and
\begin{eqnarray}\label{condition-cor}
\frac{|c_1|}{\alpha}Z^{p-1}<\inf_{u\in H^1_0(\Omega),
\|u\|_{0,\Omega}=1}\left( \int_{\Omega} |\nabla
u|^2/\int_{\Omega}|u|^{2p}\right).
\end{eqnarray}
\end{itemize}

\item[(iii)] $\mathcal{N}_1(t)\in Pol(p_1,(c_1,c_2))$ for some
$p_1\in [0,2)$ and $\mathcal{N}'_1(t)\in Pol(p_2,(c_1,c_2))$ for
some $p_2\in [0,1)$.

\item[(iv)] $K\in L^2(\tilde{\Omega})$, where $\tilde{\Omega}
=\{x-y:x,y\in\Omega\}$. Moreover, $K$ is some nonnegative even
function and $q\in[1,3/2)$.
\end{itemize}
Note that these assumptions are satisfied by typical physical models
for micro-structures of the matter (see, e.g.,
\cite{blanc-cances05,bris03,chen-zhou08,dalfovo-etal99,lieb81}) and
condition (\ref{condition-cor}) was first appeared in
\cite{blanc-cances05}.

It is known that under assumptions (i)-(iv), there exists a
nonnegative minimizer of energy functional (\ref{energy})
\cite{chen-gong-zhou09,lieb81,zhou07}.  Moreover, $E(v)$ is bounded
below over $\mathcal{A}$ under these assumptions
\cite{chen-gong-zhou09,zhou07}, namely, there exist positive
constants $C$ and $b$ such that
\begin{eqnarray}\label{bounded}
E(v)\geq C^{-1}\|\nabla v\|^2_{0,\Omega}-b, \quad\forall~ v\in
\mathcal{A}.
\end{eqnarray}
In general, however, the uniqueness of the nonnegative ground state
solution is unknown, of which the main reason is that energy
functional (\ref{energy}) is nonconvex with respect to $\rho=u^2$
for almost all molecular models of practical interest. As a result,
we introduce the set of ground state solutions by
\begin{eqnarray}\label{set-solutions}
\mathcal{U}=\{u\in\mathcal{A}:E(u)=\min_{v\in\mathcal{A}}E(v)\}.
\end{eqnarray}
If $u\in\mathcal{U}$ is a ground state solution, then there exists a
corresponding Lagrange multiplier $\lambda\in \mathbb{R}$ such that
$(\lambda,u)$ solves (\ref{weak}) and satisfies
\begin{eqnarray*}\label{eigen-energy-}
Z\lambda=E(u)+\int_{\Omega}\left(\mathcal{N}_1
(u^2(x))u^2(x)-\mathcal{E}(u^2(x))\right)dx
+(1-\frac{1}{2q})D_K(u^{2q},u^{2q}).
\end{eqnarray*}
We define the set of ground state eigenvalues by
\begin{eqnarray*}
\Lambda=\{\lambda\in\mathbb{R}:(\lambda,u)\mbox{ solves
\eqref{weak}}, ~u\in\mathcal{U}\}.
\end{eqnarray*}

The following estimate of the nonlinear term will be used in our
analysis.
\begin{lemma}
Let $\chi, w\in H^1(\Omega)$ satisfy
$\|\chi\|_{1,\Omega}+\|w\|_{1,\Omega}\leq \bar{C}$ for some constant
$\bar{C}$. Then there exists a constant ${\tilde C}>0$ depending on
${\bar C}$ such that
\begin{eqnarray}\label{n-estimate}
\int_{\Omega} \big(\mathcal{N}(\chi^2)\chi-\mathcal{N}(w^2)w\big) v
\le {\tilde C}  \|\chi-w\|_{1,\Omega}\|v\|_{1,\Omega}\quad \forall~
v\in H_0^1(\Omega).
\end{eqnarray}
\end{lemma}

\begin{proof}
We first prove that (\ref{n-estimate}) holds when $\mathcal{N}$ is
replaced by  local term $\mathcal{N}_1$. Since there exists
$\delta\in [0,1]$ such that
\begin{eqnarray*}
\mathcal{N}_1(\chi^2)\chi-\mathcal{N}_1(w^2)w =
\big(\mathcal{N}_1(\xi^2)+2\xi^2\mathcal{N}_1'(\xi^2)\big) (\chi-w)
\end{eqnarray*}
with $\xi=\chi+\delta(w-\chi)$, we have
\begin{eqnarray*}
\int_{\Omega}
\big(\mathcal{N}_1(\chi^2)\chi-\mathcal{N}_1(w^2)w\big)v
=\int_{\Omega}
\big(\mathcal{N}_1(\xi^2)+2\xi^2\mathcal{N}_1'(\xi^2)) (\chi-w)
\big)v \quad \forall~ v\in H_0^1(\Omega).
\end{eqnarray*}
From assumption (iii), we have that for all $v\in H_0^1(\Omega)$,
there hold
\begin{eqnarray*}
\int_{\Omega}\mathcal{N}_1(\xi^2)(\chi-w)v &\lesssim &
\|\xi^{2p_1}\|_{0,3/p_1,\Omega}\|\chi-w\|_{0,6/(5-2p_1),\Omega}
\|v\|_{0,6,\Omega}\nonumber\\ &\lesssim&
\|\xi\|^{2p_1}_{0,6,\Omega}\|\chi-w\|_{1,\Omega}\|v\|_{1,\Omega}
\end{eqnarray*}
and
\begin{eqnarray*}
\int_{\Omega} 2\xi^2\mathcal{N}_1'(\xi^2)(\chi-w)v &\lesssim &
\|\xi^{2p_2+2}\|_{0,3/2,\Omega}\|\chi-w\|_{0,6,\Omega}
\|v\|_{0,6,\Omega}\\&\lesssim &\|\xi\|^{2p_2+2}_{0,6,\Omega}
\|\chi-w\|_{1,\Omega}\|v\|_{1,\Omega},
\end{eqnarray*}
where the H{\" o}lder inequality and the Sobolev inequality are
used. Note that
\begin{eqnarray*}
\|\xi\|_{0,6,\Omega}\lesssim\|\chi\|_{0,6,\Omega}+\|w\|_{0,6,\Omega}
\lesssim\|\chi\|_{1,\Omega}+\|w\|_{1,\Omega}\le\bar{C},
\end{eqnarray*}
so we get
\begin{eqnarray}\label{them-ineq1}\nonumber
\int_{\Omega}\big(\mathcal{N}_1(\chi^2)\chi
-\mathcal{N}_1(w^2)w\big)v &\lesssim&
(\|\xi\|^{2p_1}_{0,6,\Omega}+\|\xi\|^{2p_2+2}_{0,6,\Omega})
\|\chi-w\|_{1,\Omega}\|v\|_{1,\Omega} \\ &\lesssim&
\|\chi-w\|_{1,\Omega}\|v\|_{1,\Omega} \quad\forall ~v\in
H_0^1(\Omega).
\end{eqnarray}

For  nonlocal term $\mathcal{N}_2$, we obtain from assumption (iv),
the Young's inequality, and the H{\" o}lder inequality that
\begin{eqnarray*}
\|K*(\chi^{2q}-w^{2q})\|_{0,\infty,\Omega}&\lesssim&
\|K\|_{0,\tilde{\Omega}}\|\chi^{2q}-w^{2q}\|_{0,\Omega}\\
&\lesssim&\|K\|_{0,\tilde{\Omega}}\|\xi^{2q-1}\|_{0,6/5,\Omega}
\|\chi-w\|_{0,6,\Omega} \\ &\lesssim& \|\chi-w\|_{1,\Omega}.
\end{eqnarray*}
Hence, for all $v\in H^1_0(\Omega)$ we have
\begin{eqnarray}\label{proof-2-1}
\nonumber \int_{\Omega} K*(\chi^{2q}-w^{2q})\chi^{2q-1}v &\lesssim&
\|K*(\chi^{2q}-w^{2q})\|_{0,\infty,\Omega}\|\chi^{2q-1}\|_{0,\Omega}
\|v\|_{0,\Omega}\\ \nonumber &\lesssim& \|\chi-w\|_{1,\Omega}
\|\chi\|^{2q-1}_{0,2(2q-1),\Omega}\|v\|_{0,\Omega}\\
&\lesssim &\|\chi-w\|_{1,\Omega}\|v\|_{1,\Omega} ,
\end{eqnarray}
where $q\in[1,3/2)$ as in assumption (iv). Similarly, for all $v\in
H^1_0(\Omega)$, there holds
\begin{eqnarray}\label{proof-2-2}
\nonumber \int_{\Omega} K*w^{2q}(\chi^{2q-1}-w^{2q-1})v &\lesssim&
\|K*w^{2q}\|_{0,\infty,\Omega}\|\chi^{2q-1}-w^{2q-1}\|_{0,\Omega}
\|v\|_{0,\Omega}\nonumber\\ &\lesssim &
\|\xi^{2q-2}\|_{0,3/(q-1),\Omega}\|\chi-w\|_{0,6/(5-2q),\Omega}
\|v\|_{1,\Omega}\nonumber\\ &\lesssim& \|\chi-w\|_{1,\Omega}
\|v\|_{1,\Omega}.
\end{eqnarray}
Taking  (\ref{proof-2-1}), (\ref{proof-2-2}), and identity
\begin{eqnarray*}
\mathcal{N}_2(\chi^2)\chi-\mathcal{N}_2(w^2)w
=K*(\chi^{2q}-w^{2q})\chi^{2q-1} + K*w^{2q}(\chi^{2q-1}-w^{2q-1})
\end{eqnarray*}
into account, we obtain
\begin{eqnarray*}
\int_{\Omega}
\big(\mathcal{N}_2(\chi^2)\chi-\mathcal{N}_2(w^2)w\big) v
\lesssim\|\chi-w\|_{1,\Omega}\|v\|_{1,\Omega}\quad\forall~ v\in
H_0^1(\Omega),
\end{eqnarray*}
which together with (\ref{them-ineq1}) leads to (\ref{n-estimate}).
This completes the proof.
\end{proof}

\section{Finite element discretizations}\label{sec-FED} Let
$\{\mathcal{T}_h\}$ be a shape regular family of nested conforming
meshes over $\Omega$ with size $h$: there exists a constant
$\gamma^*$ such that
\begin{eqnarray*}\label{eq-regularity}
\frac{h_T}{\rho_T} \leq \gamma^* \quad\quad \forall~ T\in
\mathcal{T}_h,
\end{eqnarray*}
where, for each $T\in\mathcal{T}_h$, $h_T$ is the diameter of $T$,
$\rho_T$ is the diameter of the biggest ball contained in $T$, and
$h=\max\{h_T:T\in\mathcal{T}_h\}$. Let $\mathcal{E}_h$ denote the
set of interior faces  of $\mathcal{T}_h$. And we shall also use a
slightly abused notation that  $h$ denotes the mesh size function
defined by
$$
h(x)=h_T, ~~ x\in T~~\forall~ T\in \mathcal{T}_h.$$

Let $S^h(\Omega)\subset H^1(\Omega)$ be a corresponding family of
nested finite element spaces consisting of continuous piecewise
polynomials over $\mathcal{T}_h$ of fixed degree $n\geq 1$ and
$S_0^h(\Omega)=S^h(\Omega)\cap H_0^1(\Omega)$. Set
$V_h=S^h_0(\Omega)\cap \mathcal{A}$.

Under assumptions (i)-(iv), we can obtain the existence of
nonnegative ground state solutions in $V_h$ (see, e.g.,
\cite{zhou07}). We do not have any uniqueness result for this
discrete problem since the energy functional and the admissible set
are nonconvex. We define the set of ground state solutions in $V_h$
by
\begin{eqnarray}\label{energy-dis}
\mathcal{U}_h=\{u_h\in V_h: E(u_h)=\min_{v\in V_h}E(v)\}.
\end{eqnarray}
We have from (\ref{bounded}) that $\|u_h\|_{1,\Omega}$ is uniformly
bounded:
\begin{eqnarray}\label{bounded-dis}
\sup_{h<1,u_h\in \mathcal{U}_h}\|u_h\|_{1,\Omega}\leq C
\end{eqnarray}
with some constant $C$.

It is seen that a minimizer $u_h\in \mathcal{U}_h$ solves
\begin{eqnarray}\label{weak-dis}
\left\{ \begin{array}{rcl} \alpha(\nabla u_h,\nabla v) + (Vu_h +
\mathcal{N}(u_h^2)u_h, v) & = &
\lambda_h(u_h,v) \quad \forall~v\in S_0^h(\Omega), \\[1ex]
\|u_h\|^2_{0,\Omega} &=& Z\end{array} \right.
\end{eqnarray}
with the corresponding finite element eigenvalue $\lambda_h\in
\mathbb{R}$ satisfying
\begin{eqnarray}\label{eigen-energy}
Z\lambda_h=E(u_h)+\int_{\Omega}\left(\mathcal{N}_1(u_h^2(x))
u_h^2(x)-\mathcal{E}(u_h^2(x))\right)dx
+(1-\frac{1}{2q})D_K(u_h^{2q},u_h^{2q}).
\end{eqnarray}
Define
\begin{gather*}
\Lambda_h=\{\lambda_h\in\mathbb{R}:(\lambda_h,u_h)\mbox{ solves
\eqref{weak-dis}}, u_h\in \mathcal{U}_h\}.
\end{gather*}

A priori error analysis for (\ref{weak-dis}) has been shown in
\cite{chen-gong-zhou09}. To carry out a posteriori error analysis,
we need the following result.

\begin{lemma}
There hold
\begin{eqnarray}\label{T-n-estimate}
h_T\|\mathcal{N}(u_h^2)u_h\|_{0,T}\lesssim\|u_h\|_{1,T}\quad
\forall~ T\in \mathcal{T}_h.
\end{eqnarray}
\end{lemma}

\begin{proof}
It is obvious that
\begin{eqnarray*}
h_T\|\mathcal{N}_1(u_h^2)u_h\|_{0,T} \lesssim\|u_h\|_{1,T}
\end{eqnarray*}
holds for $p_1=0$ in assumption (iii). By the H\"{o}lder inequality
and the inverse inequality,
\begin{eqnarray*}
h_T\|\mathcal{N}_1(u_h^2)u_h\|_{0,T} &\lesssim &
h_T\|\mathcal{N}_1(u_h^2)\|_{0,3,T}\|u_h\|_{0,6,T}\lesssim
h_T\|u_h\|_{0,6p_1,T}^{2p_1}\|u_h\|_{0,6,T}\\ &\lesssim &
h_T^{2-p_1}\|u_h\|^{2p_1}_{1,T}\|u_h\|_{1,T}\lesssim\|u_h\|_{1,T}
\end{eqnarray*}
for $p_1\in[1,2)$ and
\begin{eqnarray*}
h_T\|\mathcal{N}_1(u_h^2)u_h\|_{0,T}\lesssim h_T\|u_h^{2p_1}
\|_{0,3/p_1,T}\|u_h\|_{1,6/(3-2p_1),T} \lesssim \|u_h\|_{1,T}
\end{eqnarray*}
for $p_1\in (0,1)$. Combining with the estimate of $\mathcal{N}_2$
as follows
\begin{eqnarray*}
h_T\|\mathcal{N}_2(u_h^2)u_h\|_{0,T} &\lesssim &
h_T\|K*u_h^{2q}\|_{0,\infty,T}\|u_h^{2q-1}\|_{0,T}\\ &\lesssim&
h_T\|K\|_{0,\tilde{\Omega}}\|u_h^{2q}\|_{0,\Omega}
\|u_h\|^{2q-1}_{1,T} \lesssim\|u_h\|_{1,T},
\end{eqnarray*}
where assumption (iv) is used, we obtain (\ref{T-n-estimate}). This
completes the proof.
\end{proof}

Let $\mathbb{T}$ denote the class of all conforming refinements by
the bisection of an initial triangulation $\mathcal{T}_0$. For
$\mathcal{T}_h \in \mathbb{T}$ and any $u_h\in \mathcal{U}_h$ we
define  element residual $\mathcal{R}_T(u_h)$ and  jump residual
$J_e(u_h)$ by
\begin{gather*}
\mathcal{R}_T(u_h) = \lambda_h u_h +\alpha \Delta
u_h-Vu_h-\mathcal{N}(u_h^2)u_h \quad \mbox{in}~
T\in \mathcal{T}_h,\\[1ex]
J_e(u_h) =\alpha \nabla u_h |_{T_1}\cdot\overrightarrow{n_1} +
\alpha \nabla u_h |_{T_2}\cdot\overrightarrow{n_2}=[[\alpha \nabla
u_h]]_e \cdot \overrightarrow{n_1} \quad \mbox{on}~ e\in
\mathcal{E}_h,
\end{gather*}
where $T_1$ and $T_2$ are elements in $\mathcal{T}_h$ which share
$e$ and $\overrightarrow{n_i}$ is the outward normal vector of $T_i$
on $E$ for $i=1,2$. Let $\omega_h(e)$ be the union of elements which
share $e$ and $\omega_h(T)$ be the union of elements sharing a side
with $T$.

For $T\in \mathcal{T}_h$, we define  local error indicator
$\eta_h(u_h, T)$ by
\begin{eqnarray}\label{Gerror-indicator}
\eta^2_h(u_h, T) = h_T^2\|\mathcal{R}_{T}(u_h)\|_{0,T}^2 +
\sum_{e\in \mathcal{E}_h,e\subset\partial T} h_e
\|J_e(u_h)\|_{0,e}^2.
\end{eqnarray}
Given a subset $\omega \subset \Omega$, we define error estimator
$\eta_h(u_h, \omega)$ by
\begin{eqnarray*}\label{Gerror-estimator}
\eta^2_h(u_h, \omega) = \sum_{T\in \mathcal{T}_h, T\subset \omega}
\eta^2_h(u_h, T).
\end{eqnarray*}

The following result will be used in our convergence analysis though
it looks rough.
\begin{proposition}\label{estimator-stability}
Let $\mathcal{T}_h\in \mathbb{T}$. If $(\lambda_h, u_h)$ is a
solution of (\ref{weak-dis}), then
\begin{eqnarray*}\label{eta-stability}
\eta_h(u_h,T)\lesssim \|u_h\|_{1,\omega_h(T)} \quad\forall~ T\in
\mathcal{T}_h
\end{eqnarray*}
and
\begin{eqnarray*}\label{eta-stability-}
\eta_h(u_h,\Omega)\leq C_{\eta},
\end{eqnarray*}
where the uniform constant $C_{\eta}>0$ depends only on the data and
the mesh regularity.
\end{proposition}
\begin{proof}
We first analyze the element residual. Note that
\begin{eqnarray*}
h_T\|\mathcal{R}_T(u_h)\|_{0,T} &=& h_T\|\lambda_hu_h+\alpha \Delta
u_h-Vu_h-\mathcal{N}(u_h^2)u_h\|_{0,T}\nonumber\\
&\lesssim &h_T\|u_h\|_{0,T}+h_T\|\Delta
u_h\|_{0,T}+h_T\|Vu_h\|_{0,T} +h_T\|\mathcal{N}(u_h^2)u_h\|_{0,T}.
\end{eqnarray*}
Using  the inverse inequality, assumption (i) and
(\ref{T-n-estimate}), we have
\begin{eqnarray*}
h_T\|\mathcal{R}_T(u_h)\|_{0,T}\lesssim \|u_h\|_{1,T},
\end{eqnarray*}
to which similar estimates are true when $T$ is replaced by any
$T'\in \omega_h(T)$.

For the jump residual, by the definition of $J_e(u_h)$ and the trace
inequality,
\begin{eqnarray*}
h_e^{1/2}\|J_e(u_h)\|_{0,e} &=& h_e^{1/2}\|\alpha \nabla u_h
\big|_{T_1}\cdot\overrightarrow{n_1} + \alpha \nabla
u_h\big|_{T_2}\cdot\overrightarrow{n_2}\|_{0,e}\nonumber\\
&\lesssim& h_e^{1/2}(\|\nabla u_h|_{T_1}\|_{0,e} +\|\nabla
u_h|_{T_2}\|_{0,e})\\ &\lesssim & \|u_h\|_{1,\omega_h(T)}.
\end{eqnarray*}
Hence
\begin{eqnarray*}
\sum_{e\in \mathcal{E}_h, e\subset
\partial T}h_e^{1/2}\|J_e(u_h)\|_{0,e}\lesssim \|u_h\|_{1,\omega_h(T)}.
\end{eqnarray*}
From \eqref{bounded-dis}, the definition of $\eta_h(u_h,T)$ and
$\eta_h(u_h,\Omega)$, we get the desired result.
\end{proof}

To present upper and lower error bounds, we introduce an oscillation
$osc_h(u_h,T)$ for any $T\in\mathcal{T}_h$ by
\begin{eqnarray*}\label{local-oscillation}
{\rm osc}^2_h(u_h,T) =
h_T^2\|\mathcal{R}_{T}(u_h)-\overline{\mathcal{R}_{T}(u_h)}\|_{0,T}^2,
\end{eqnarray*}
where $\overline{\mathcal{R}_{T}(u_h)}\in P_{n-1}$ denotes the $L^2$
projection of $\mathcal{R}_{T}(u_h)$. For a subset $\omega \subset
\Omega$, we define
\begin{eqnarray*}\label{Goscilliation}
{\rm osc}^2_h(u_h, \omega) = \sum_{T\in \mathcal{T}_h, T \subset
\omega} {\rm osc}^2_h(u_h, T).
\end{eqnarray*}

We have have a standard argument that (see Appendix for a  proof)
\begin{theorem}\label{theorem-bound}
Let $(\lambda,u)$ be a regular ground state solution of
\eqref{weak}. If $(\lambda_h,u_h)$ is sufficiently close to
$(\lambda,u)$, then
\begin{eqnarray*}\label{eq-bound}
\eta_h(u_h,\Omega)-{\rm osc}_h(u_h,\Omega) \lesssim
|\lambda-\lambda_h| + \|u-u_h\|_{1,\Omega} \lesssim
\eta_h(u_h,\Omega) + {\rm osc}_h(u_h,\Omega).
\end{eqnarray*}
\end{theorem}\vskip 0.2cm
\begin{remark}
This result provides the standard upper and lower bounds of the
error with respect to the error estimator. However, the hypothesis
that $(\lambda,u)$ is a regular solution is somehow strong, which
can not be proved for most of the problems of practical interest
(c.f., e.g., Appendix). Anyway, it will not be used in our
convergence analysis.
\end{remark}

We define  global residual  $\mathbf{R}_h(u_h)\in H^{-1}(\Omega)$ as
follows
\begin{eqnarray*}
\langle\mathbf{R}_h(u_h),v\rangle = \lambda_h(u_h,v)-(\alpha \nabla
u_h,\nabla v)-(V u_h,v)-(\mathcal{N}(u_h^2)u_h,v) \quad \forall~
v\in H^1_0(\Omega)
\end{eqnarray*}
and see that
\begin{eqnarray*}\label{residual-jump}
\langle\mathbf{R}_h(u_h),v\rangle=\sum_{T\in \mathcal{T}_h}\left(
\int_{T}\mathcal{R}_T(u_h)v-\sum_{e\in \mathcal{E}_h,
e\subset\partial T}\int_{e}J_e(u_h)v\right)\quad \forall~ v\in
H^1_0(\Omega).
\end{eqnarray*}

The global residual  can be estimated by the local error indicators
in the following sense.
\begin{theorem}\label{upp-bou-resi} If $(\lambda_h,
u_h)\in \mathbb{R}\times V_h$ is a solution of (\ref{weak-dis}),
then
\begin{eqnarray*}\label{upper-residual}
|\langle\mathbf{R}_h(u_h),v\rangle| \lesssim \sum_{T\in
\mathcal{T}_h}\eta_h(u_h,T)\|v\|_{1,\omega_h(T)}~~\forall ~v\in
H^1_0(\Omega).
\end{eqnarray*}
\end{theorem}
\begin{proof}
Let $v\in H^1_0(\Omega)$ and $v_h\in S^h_0(\Omega)$ be the
Cl\'{e}ment interpolant of $v$ satisfying
\begin{eqnarray*}
\|v-v_h\|_{0,T}\lesssim h_T\|\nabla v\|_{0,\omega_h(T)}\quad
\textnormal{and} \quad \|v-v_h\|_{0,\partial T}\lesssim
h_T^{1/2}\|\nabla v\|_{0,\omega_h(T)}.
\end{eqnarray*}
Due to $\langle\mathbf{R}_h(u_h),v_h\rangle=0$, we obtain
\begin{eqnarray*}
&&|\langle\mathbf{R}_h(u_h),v\rangle|
=|\langle\mathbf{R}_h(u_h),v-v_h\rangle|\nonumber\\
&\leq& \sum_{T\in \mathcal{T}_h}
\left(\|\mathcal{R}_T(u_h)\|_{0,T}\|v-v_h\|_{0,T}+\sum_{e\in
\mathcal{E}_h,e\subset\partial T}\|J_e(u_h)\|_{0,e}
\|v-v_h\|_{0,e}\right)\nonumber\\
&\lesssim&\sum_{T\in \mathcal{T}_h}\left(\|h\mathcal{R}_T(u_h)
\|_{0,T}\|v\|_{1,\omega_h(T)}+\sum_{e\in
\mathcal{E}_h,e\subset\partial T}\|h^{1/2}J_e(u_h)
\|_{0,e}\|v\|_{1,\omega_h(e)}\right) \nonumber\\ &\lesssim&
\sum_{T\in\mathcal{T}_h}\eta_h(u_h,T)\|v\|_{1,\omega_h(T)}.
\end{eqnarray*}
This completes the proof.
\end{proof}

\section{Convergence of adaptive finite element computations}
\label{sec-convergence} \setcounter{equation}{0} We shall first
recall the adaptive finite element algorithm.  For convenience, we
shall replace the subscript $h$ (or $h_k$) by an iteration counter
$k$ of the adaptive algorithm afterwards. Given an initial
triangulation $\mathcal{T}_0$, we can generate a sequence of nested
conforming triangulations $\mathcal{T}_k$ using the following loop:
$$
\mbox{Solve}~\rightarrow~\mbox{Estimate}~\rightarrow~
\mbox{Mark}~\rightarrow~\mbox{Refine}.
$$
More precisely, to get $\mathcal{T}_{k+1}$ from $\mathcal{T}_k$ we
first solve the discrete equation to get $\mathcal{U}_k$ on
$\mathcal{T}_k$. The error is estimated by any $u_k\in
\mathcal{U}_k$ and used to mark a set of elements that are to be
refined. Elements are refined in such a way that the triangulation
is still shape regular and conforming.

Here, we shall not discuss the step ``Solve", which deserves a
separate investigation. We assume that  solutions of finite
dimensional problems can be solved to any accuracy efficiently. The
procedure ``Estimate" determines the element indicators for all
elements $T\in \mathcal{T}_k$. A posteriori error estimators are an
essential part of this step, which have been investigated in the
previous section.  In the following discussion, we use
$\eta_k(u_k,T)$ defined by (\ref{Gerror-indicator}) as the a
posteriori error estimator. Depending on the relative size of the
element indicators, these quantities are later used by the procedure
``Mark" to mark elements in $\mathcal{T}_k$ and thereby create a
subset of elements to be refined. The only requirement we make on
this step is that the set of marked elements $\mathcal{M}_k$
contains at least one element of $\mathcal{T}_k$ holding the largest
value estimator \cite{garau-morin-zuppa09,garau-morin09}. Namely,
there exists one element $T^{\max}_k\in \mathcal{M}_k$ such that
\begin{eqnarray}\label{mark-pri}
\eta_k(u_k,T^{\max}_k)=\max_{T\in \mathcal{T}_k}\eta_k(u_k,T).
\end{eqnarray}
It is easy to check that the most commonly used marking strategies,
e.g., Maximum strategy, Equidistribution strategy, and D\"{o}rfler's
strategy fulfill this condition. Finally, the marked elements are
refined to force the error reduction by the procedure ``Refine". The
basic algorithm in this step is the tetrahedral bisection, with the
data structure named marked tetrahedron, the tetrahedra are
classified into 5 types and the selection of refinement edges
depends only on the type and the ordering of vertices for the
tetrahedra \cite{arnold-mukherjee-pouly00}. Note that a few more
elements $T\in\mathcal{T}_k ~\backslash ~\mathcal{M}_k$ are
partitioned to maintain mesh conformity. It is worth mentioning that
we do not assume to enforce the so-called interior node property.

The adaptive finite element algorithm without oscillation marking is
stated as follows: \vskip 0.1cm

\begin{algorithm}\label{algorithm-AFEM}~~

\begin{enumerate}
\item Pick  any initial mesh $\mathcal{T}_0$, and let $k=0$.
\item Solve the system on $\mathcal{T}_k$ to get discrete
solutions $\mathcal{U}_k$.
\item Choose any $u_k\in \mathcal{U}_k$ and compute local error
indictors $\eta_k(u_k,T)~~\forall~ T\in \mathcal{T}_k$.
\item Construct $\mathcal{M}_k \subset \mathcal{T}_k$ by a marking
strategy that satisfies (\ref{mark-pri}).
\item Refine $\mathcal{T}_k$ to get a new conforming mesh
$\mathcal{T}_{k+1}$.
\item Let $k=k+1$ and go to 2.
\end{enumerate}
\end{algorithm}

The purpose of this paper is to prove that Algorithm
\ref{algorithm-AFEM} generates a sequence of adaptive finite element
solutions  which converge to some ground state solutions of
(\ref{set-solutions}). More precisely, we shall prove that
\begin{gather*}
\lim_{k\rightarrow\infty}
\textnormal{dist}_{H^1}(\mathcal{U}_k,\mathcal{U})=0,\\
\lim_{k\rightarrow\infty}\textnormal{dist}(\Lambda_k,\Lambda)=0,
\end{gather*}
 where
\begin{eqnarray*}
\textnormal{dist}_{H^1}(F,G) = \sup_{f\in F}\inf_{g\in
G}\|f-g\|_{1,\Omega}
\end{eqnarray*}
for any $F,G\subset H^1(\Omega)$, and
\begin{eqnarray*}
\textnormal{dist}(A,B) =\sup_{a\in A}\inf_{b\in B}|a-b|
\end{eqnarray*}
for any $A,B\subset \mathbb{R}$.

We first show that the adaptive finite element approximations are
convergent. Given an initial mesh $\mathcal{T}_0$, Algorithm
\ref{algorithm-AFEM} generates a sequence of meshes
$\mathcal{T}_1,\mathcal{T}_2,\cdots$, and associated discrete
subspaces
\begin{eqnarray*}
S_0^{h_0}(\Omega)\subsetneq S_0^{h_1}(\Omega)\subsetneq \cdots
\subsetneq S_0^{h_n}(\Omega)\subsetneq S_0^{h_{n+1}}(\Omega)
\subsetneq \cdots \subsetneq S_{\infty}(\Omega) \subseteq
H_0^1(\Omega),
\end{eqnarray*}
where $\displaystyle S_{\infty}(\Omega)=\overline{\cup
S_0^{h_k}(\Omega)}^{H_0^1(\Omega)}$. It is obvious that
$S_{\infty}(\Omega)$ is a Hilbert space with the inner product
inherited from $H_0^1(\Omega)$, and there holds
\begin{eqnarray}\label{space-appro}
\lim_{k\rightarrow\infty}\inf_{v_k\in S_0^{h_k}(\Omega)}
\|v_k-v_{\infty}\|_{1,\Omega}=0 \quad\forall~ v_{\infty}\in
S_{\infty}(\Omega).
\end{eqnarray}
We set $V_{\infty}=S_{\infty}(\Omega)\cap \mathcal{A}$.

Under assumptions (i)-(iv), the existence of minimizers of energy
functional \eqref{energy} in $V_{\infty}$ can be obtained. Similar
to (\ref{set-solutions}) and (\ref{energy-dis}), we introduce the
set of minimizers by
\begin{eqnarray*}\label{set-solutions-infty}
\mathcal{U_{\infty}}=\{u\in V_{\infty}:E(u)=\min_{v\in V_{\infty}}
E(v)\}.
\end{eqnarray*}
We see that $u_{\infty}\in\mathcal{U}_{\infty}$ solves
\begin{eqnarray}\label{weak-dis-infty}
\left\{ \begin{array}{rcl} \alpha(\nabla u_{\infty},\nabla v) +
(Vu_{\infty} + \mathcal{N}(u_{\infty}^2)u_{\infty}, v) & = &
\lambda_{\infty}(u_{\infty},v) \quad \forall~v\in
S_{\infty}(\Omega),
\\[1ex] \|u_{\infty}\|^2_{0,\Omega} &=& Z\end{array} \right.
\end{eqnarray}
with the corresponding eigenvalue $\lambda_{\infty}\in \mathbb{R}$
satisfying
\begin{eqnarray}\label{eigen-energy-infty}
Z\lambda_{\infty}=E(u_{\infty})+\int_{\Omega}\left(\mathcal{N}_1
(u_{\infty}^2(x))u_{\infty}^2(x)-\mathcal{E}(u_{\infty}^2(x))\right)dx
+(1-\frac{1}{2q})D_K(u_{\infty}^{2q},u_{\infty}^{2q}),
\end{eqnarray}
and we define
\begin{gather*}
\Lambda_{\infty}=\{\lambda_{\infty}\in\mathbb{R}:(\lambda_{\infty},
u_{\infty})\mbox{ solves \eqref{weak-dis-infty}},
~u_{\infty}\in\mathcal{U}_{\infty}\}.
\end{gather*}

\begin{theorem}\label{theorem-con-infty}
If $\{\mathcal{U}_k\}_{k\in \mathbb{N}}$ is the sequence of adaptive
finite element approximations generated by Algorithm
\ref{algorithm-AFEM}, then
\begin{gather*}\label{con-E-infty}
\lim_{k\rightarrow\infty}E_k = \min_{v\in V_{\infty}} E(v),\\[1ex]
\label{con-u-infty} \lim_{k\rightarrow\infty}
\textnormal{dist}_{H^1}(\mathcal{U}_k,\mathcal{U}_{\infty})=0,
\end{gather*}
where $E_k=E(v) (v\in\mathcal{U}_k)$. Moreover, there holds
\begin{eqnarray*}\label{con-eigenvalue-infty}
\lim_{k\rightarrow\infty}\textnormal{dist}(\Lambda_k,\Lambda_{\infty})=0.
\end{eqnarray*}
\end{theorem}

\begin{proof}
Following \cite{zhou04,zhou07} (see also \cite{chen-gong-zhou09}),
let $u_k\in\mathcal{U}_k$  be such that $(\lambda_{k},u_{k})$ solves
(\ref{weak-dis}) in $\mathbb{R}\times V_{k}$ for $k=1,2,\cdots$, and
$\{u_{k_m}\}_{m\in \mathbb{N}}$ be any subsequence of $\{u_k\}_{k\in
\mathbb{N}}$ with $1\leq k_1<k_2<\cdots <k_m<\cdots$.

Note that (\ref{bounded-dis}) and the Banach-Alaoglu Theorem yield
that there exist a weakly convergent subsequence
$\{u_{k_{m_j}}\}_{j\in \mathbb{N}}$ and $u_{\infty}\in
S_{\infty}(\Omega)$ satisfying
\begin{eqnarray}\label{con-weak-proof}
u_{k_{m_j}}\rightharpoonup u_{\infty} \quad
\mbox{in}~~H_0^1(\Omega),
\end{eqnarray}
we need only to prove
\begin{eqnarray}\label{temp1}
E(u_{\infty})=\min_{v\in V_{\infty}}E(v),
\end{eqnarray}
\begin{eqnarray}\label{temp2}
\lim_{j\to\infty}\|u_{k_{m_j}}-u_{\infty}\|_{1,\Omega}=0,
\end{eqnarray}
and
\begin{eqnarray}\label{temp3}
\lim_{j\to\infty}|\lambda_{k_{m_j}}-\lambda_{\infty}|=0,
\end{eqnarray}
where $(\lambda_{k_{m_j}},u_{k_{m_j}})$ solves \eqref{weak-dis} and
$(\lambda_{\infty},u_{\infty})$ solves \eqref{weak-dis-infty}.

Since $H_0^1(\Omega)$ is compactly imbedded in $L^p(\Omega)$ for
 $p\in[2,6)$, by passing to a further
subsequence, we may assume that $u_{k_{m_j}}\to u_{\infty}$ strongly
in $L^p(\Omega)$ as $j\rightarrow \infty$. Thus we can derive
\begin{gather*}
\lim_{j\to\infty}\int_{\Omega}\mathcal{E}(u_{k_{m_j}}^2)
=\int_{\Omega}\mathcal{E}(u_{\infty}^2), \\[1ex]
\lim_{j\to\infty}D_K(u^{2q}_{k_{m_j}},
u^{2q}_{k_{m_j}})=D_K(u_{\infty}^{2q},u_{\infty}^{2q}),
\end{gather*}
and hence
\begin{eqnarray}\label{semi-lower}
\liminf_{j\to\infty}E(u_{k_{m_j}})\geq E(u_{\infty}).
\end{eqnarray}

Note that \eqref{space-appro} implies that $\{u_{k_{m_j}}\}$ is a
minimizing sequence for the energy functional in
$S_{\infty}(\Omega)$, which together with (\ref{semi-lower}) and the
fact that $\{u_{k_{m_j}}\}$ converge to $u_{\infty}$ strongly in
$L^2(\Omega)$ leads to $u_{\infty}\in\mathcal{U}_{\infty}$, namely,
\begin{eqnarray*}\label{engery-subconv}
\lim_{j\to\infty}E(u_{k_{m_j}})=E(u_{\infty})=\min_{v\in
V_{\infty}}E(v).
\end{eqnarray*}
Consequently, we obtain that each term of $E(v)$ converges and in
particular
\begin{eqnarray*}
\lim_{j\to\infty}\|\nabla u_{k_{m_j}}\|_{0,\Omega} = \|\nabla
u_{\infty}\|_{0,\Omega}.
\end{eqnarray*}
Using (\ref{con-weak-proof}) and the fact that $H_0^1(\Omega)$ is a
Hilbert space under the norm $\|\nabla\cdot\|_{0,\Omega}$, we have
\begin{eqnarray*}
\lim_{j\to\infty}\|\nabla(u_{k_{m_j}}-u_{\infty})\|_{0,\Omega}=0,
\end{eqnarray*}
which implies \eqref{temp2}

Using  (\ref{eigen-energy}), (\ref{eigen-energy-infty}),
(\ref{temp1}) and (\ref{temp2}), we immediately obtain
(\ref{temp3}). This completes the proof.
\end{proof}

Following the ideas in
\cite{garau-morin-zuppa09,garau-morin09,morin-siebrt-veeser08,siebert08},
we then prove the convergence of the a posteriori error estimators
and the weak convergence of residual $\mathbf{R}_k(u_{k})$, which
will be used to prove that the adaptive finite element
approximations converge to the  ground state solutions. Given the
sequence $\{\mathcal{T}_k\}_{k\in\mathbb{N}}$, for each $k\in
\mathbb{N}$ we define
\begin{eqnarray*}
\mathcal{T}^{+}_k=\{T\in \mathcal{T}_k:T\in \mathcal{T}_l,~ \forall~
l\geq k\}\quad \textnormal{and} \quad \mathcal{T}^{0}_k=
\mathcal{T}_k\setminus\mathcal{T}^{+}_k.
\end{eqnarray*}
Namely, $\mathcal{T}^+_k$ is the set of elements of $\mathcal{T}_k$
that are not refined and $\mathcal{T}^0_k$ consists of those
elements which will eventually be refined. Set
\begin{eqnarray*}
\Omega^+_k=\cup_{T\in \mathcal{T}^+_k}\omega_k(T)\quad
\textnormal{and} \quad \Omega^0_k=\cup_{T\in
\mathcal{T}^0_k}\omega_k(T).
\end{eqnarray*}

Note that the  mesh size function $h_k\equiv h_k(x)$ associated with
$\mathcal{T}_k$ is monotonically decreasing and bounded from below
by 0, we have that
\begin{eqnarray*}
h_{\infty}(x)=\lim_{k\rightarrow \infty}h_k(x)
\end{eqnarray*}
is well-defined for almost all $x\in \Omega$ and hence defines a
function in $L^{\infty}(\Omega)$. Moreover, the convergence is
uniform \cite{morin-siebrt-veeser08}.

\begin{lemma}\label{domain-limit}
If $\{h_k\}_{k\in\mathbb{N}}$ is the sequence of mesh size functions
generated by Algorithm \ref{algorithm-AFEM}, then
\begin{eqnarray*}
\lim_{k\rightarrow \infty}\|h_k-h_{\infty}\|_{0,\infty,\Omega}=0
\end{eqnarray*}
and
\begin{eqnarray*}
\lim_{k\rightarrow
\infty}\|h_k\chi_{\Omega^0_k}\|_{0,\infty,\Omega}=0,
\end{eqnarray*}
where  $\chi_{\Omega^0_k}$ is the characteristic function of
$\Omega^0_k$
\end{lemma}

\begin{lemma}\label{marking-limit}
If $\{u_k\}_{k\in \mathbb{N}}$ is a sequence of ground state
solutions of (\ref{weak-dis}) in $\{V_k\}_{k\in \mathbb{N}}$
obtained by Algorithm \ref{algorithm-AFEM}, then
\begin{eqnarray*}
\lim_{k\rightarrow \infty}\max_{T\in \mathcal{M}_{k}}
\eta_{k}(u_{k},T)=0.
\end{eqnarray*}
\end{lemma}

\begin{proof}
We see from the proof of Theorem \ref{theorem-con-infty} that for
any subsequence $\{u_{k_m}\}$ of $\{u_k\}$, there exist a convergent
subsequence $\{u_{k_{m_j}}\}$ and
$u_{\infty}\in\mathcal{U}_{\infty}$ such that
\begin{eqnarray*}
u_{k_{m_j}} \rightarrow u_{\infty}\quad\mbox{in}~~H_0^1(\Omega).
\end{eqnarray*}
Now it is only necessary for us to prove that
\begin{eqnarray*}
\lim_{j\rightarrow \infty}\max_{T\in \mathcal{M}_{k_{m_j}}}
\eta_{k_{m_j}}(u_{k_{m_j}},T)=0.
\end{eqnarray*}
In order not to clutter the notation, we shall denote by
$\{u_k\}_{k\in\mathbb{N}}$ the subsequence
$\{u_{k_{m_j}}\}_{m\in\mathbb{N}}$, and by
$\{\mathcal{T}_k\}_{k\in\mathbb{N}}$ the sequence
$\{\mathcal{T}_{k_{m_j}}\}_{m\in\mathbb{N}}$.

Let $T_k\in \mathcal{M}_k$ be such that
\begin{eqnarray*}
\eta_k(u_k,T_k)=\max_{T\in \mathcal{M}_k}\eta_k(u_k,T).
\end{eqnarray*}
Using Proposition \ref{estimator-stability}, we obtain
\begin{eqnarray}\label{eta-ineq-1}
\eta_k(u_k,T_k)\lesssim \|u_k\|_{1,\omega_k(T_k)}\leq
\|u_k-u_{\infty}\|_{1,\Omega}+\|u_{\infty}\|_{1,\omega_k(T_k)}.
\end{eqnarray}
Since $T_k\in \mathcal{M}_k\subset \mathcal{T}_k^0$, we have
\begin{eqnarray*}
|\omega_k(T_k)|\lesssim
h_{T_k}^3\leq\|h_k\chi_{\Omega^0_k}\|^3_{0,\infty,\Omega}\rightarrow
0 \quad\textnormal{as}~~ k\rightarrow\infty,
\end{eqnarray*}
where Lemma \ref{domain-limit} is used. From Theorem
\ref{theorem-con-infty}, we have that the first term in the right
hand side of (\ref{eta-ineq-1}) tends to zero, too. This completes
the proof.
\end{proof}

\begin{lemma}\label{R-limit}
If $\{u_k\}_{k\in \mathbb{N}}$ is a sequence of ground state
solutions of (\ref{weak-dis}) in $\{V_k\}_{k\in \mathbb{N}}$
obtained by Algorithm \ref{algorithm-AFEM}, then
\begin{eqnarray*}
\lim_{k\rightarrow \infty}\langle
\mathbf{R}_{k}(u_{k}),v\rangle=0\quad\forall~ v\in H^1_0(\Omega).
\end{eqnarray*}
\end{lemma}

\begin{proof}
Using similar arguments as that in the proof of Theorem
\ref{theorem-con-infty},  for any subsequence $\{u_{k_m}\}$ of
$\{u_k\}$, there exist a convergence subsequence $\{u_{k_{m_j}}\}$
and $u_{\infty}\in\mathcal{U}_{\infty}$ such that
\begin{eqnarray*}
u_{k_{m_j}} \rightarrow u_{\infty}\quad\mbox{in}~~H_0^1(\Omega),
\end{eqnarray*}
and we  need only to prove
\begin{eqnarray*}
\lim_{j\rightarrow \infty}\langle
\mathbf{R}_{k_{m_j}}(u_{k_{m_j}}),v\rangle=0\quad\forall~ v\in
H^1_0(\Omega).
\end{eqnarray*}
Since $H^2_0(\Omega)$ is dense in $H^1_0(\Omega)$,  it is sufficient
to prove
\begin{eqnarray}\label{R-limit-eq2}
\lim_{j\rightarrow \infty}\langle
\mathbf{R}_{k_{m_j}}(u_{k_{m_j}}),v\rangle=0 \quad\forall ~v\in
H^2_0(\Omega).
\end{eqnarray}
For simplicity of notation, we denote by $\{u_k\}_{k\in\mathbb{N}}$
the subsequence $\{u_{k_m}\}_{m\in\mathbb{N}}$, and by
$\{\mathcal{T}_k\}_{k\in\mathbb{N}}$ the sequence
$\{\mathcal{T}_{k_m}\}_{m\in\mathbb{N}}$.

Let $v_k\in V_k$ be the Lagrange's interpolation of $v$. Since
\begin{eqnarray*}
\langle \mathbf{R}_{k}(u_{k}),v_k\rangle=0,
\end{eqnarray*}
 we have from Theorem \ref{upp-bou-resi} that
\begin{eqnarray*}
&&|\langle \mathbf{R}_{k}(u_{k}),v\rangle|=|\langle
\mathbf{R}_{k}(u_{k}),v-v_k\rangle|\leq \sum_{T\in
\mathcal{T}_k}\eta_k(u_k,T)\|v-v_k\|_{1,\omega_k(T)}.
\end{eqnarray*}
Let $n\in \mathbb{N}$ and $k>n$. By definition,
$\mathcal{T}_n^{+}\subset \mathcal{T}_k^{+}\subset \mathcal{T}_k$.
Thus we have
\begin{eqnarray*}
|\langle \mathbf{R}_k(u_{k}),v\rangle|&\leq&\sum_{T\in
\mathcal{T}_n^{+}}\eta_k(u_k,T)\|v-v_k\|_{1,\omega_k(T)}
+\sum_{T\in\mathcal{T}_k\setminus \mathcal{T}_n^{+}}\eta_k(u_k,T)
\|v-v_k\|_{1,\omega_k(T)}\nonumber\\
&\leq&\eta_k(u_k,\mathcal{T}_n^{+})\|v-v_k\|_{1,\Omega_n^{+}(T)}
+\eta_k(u_k,\mathcal{T}_k\setminus
\mathcal{T}_n^{+})\|v-v_k\|_{1,\Omega_n^{0}(T)}.
\end{eqnarray*}
Using Proposition \ref{estimator-stability}, we get
\begin{eqnarray*}
\eta_k(u_k,\mathcal{T}_k\setminus \mathcal{T}_n^{+})\leq
\eta_k(u_k,\mathcal{T}_k)\leq C_{\eta},
\end{eqnarray*}
which together with the interpolation estimate yields
\begin{eqnarray}\label{R-limit-eq3}
|\langle
\mathbf{R}_{k}(u_{k}),v\rangle|\lesssim(\eta_k(u_k,\mathcal{T}_n^{+})
+C_{\eta}\|h_n\chi_{\Omega^0_n}\|_{0,\infty,\Omega})
\|v\|_{2,\Omega}.
\end{eqnarray}

Now we shall use (\ref{R-limit-eq3}) to prove (\ref{R-limit-eq2}).
Let $\varepsilon>0$ be arbitrary. Lemma \ref{domain-limit} implies
that there exists $n\in \mathbb{N}$ such that
\begin{eqnarray}\label{R-limit-eq4}
C_{\eta}\|h_n\chi_{\Omega^0_n}\|_{0,\infty,\Omega}<\varepsilon.
\end{eqnarray}
Since $\mathcal{T}_n^{+}\subset \mathcal{T}_k^{+}\subset
\mathcal{T}_k$ and the marking strategy (\ref{mark-pri}) is
reasonable, we arrive at
\begin{eqnarray*}
\eta_k(u_k,\mathcal{T}_n^{+})\leq (\# \mathcal{T}_n^{+})^{1/2}
\max_{T\in \mathcal{T}_n^{+}}\eta_k(u_k,T)\leq
(\#\mathcal{T}_n^{+})^{1/2}\max_{T\in \mathcal{M}_k}\eta_k(u_k,T).
\end{eqnarray*}
By Lemma \ref{marking-limit}, we can select $N\geq n$ such that
\begin{eqnarray}\label{R-limit-eq5}
\eta_k(u_k,\mathcal{T}_n^{+})<\varepsilon ~~\forall~ k>N.
\end{eqnarray}
Thus we obtain (\ref{R-limit-eq2}) by combining (\ref{R-limit-eq3}),
(\ref{R-limit-eq4}) and (\ref{R-limit-eq5}). This completes the
proof.
\end{proof}

Finally, we  prove the main result of this paper.

\begin{theorem}\label{theorem-convergence}
Given a sufficiently fine initial mesh $\mathcal{T}_0$. If
$\{\mathcal{U}_k\}_{k\in\mathbb{N}}$ is the sequence of adaptive
finite element approximations generated by Algorithm
\ref{algorithm-AFEM}, then
\begin{gather}\label{con-E}
\lim_{k\rightarrow\infty}E_k=\min_{v\in\mathcal{A}}E(v),
\\[1ex] \label{con-u}
\lim_{k\rightarrow\infty}\textnormal{dist}_{H^1}(\mathcal{U}_k,
\mathcal{U})=0,
\end{gather}
and
\begin{eqnarray}\label{con-eigenvalue}
\lim_{k\rightarrow\infty}\textnormal{dist}(\Lambda_k,\Lambda)=0.
\end{eqnarray}
\end{theorem}

\begin{proof}
We shall use a similar argument as that in the proof of Theorem
\ref{theorem-con-infty}. Let $u_k\in\mathcal{U}_k$ be such that
$(\lambda_k,u_k)$ solves \eqref{weak-dis} for $k=1,2,\cdots$. For
any subsequence $\{u_{k_m}\}$ of $\{u_k\}$, there exist a convergent
subsequence $\{u_{k_{m_j}}\}$ and
$u_{\infty}\in\mathcal{U}_{\infty}$ such that
\begin{gather*}
u_{k_{m_j}} \rightarrow u_{\infty}\quad\mbox{in}~~H_0^1(\Omega),\\[1ex]
\lambda_{k_{m_j}}\rightarrow\lambda_{\infty},
\end{gather*}
where $(\lambda_{\infty},u_{\infty})$ solves $\eqref{weak-dis}$. It
is only necessary for us to prove that $u_{\infty}\in\mathcal{U}$,
which derives \eqref{con-E} and \eqref{con-u} directly, and implies
\eqref{con-eigenvalue} by noting \eqref{weak} and \eqref{weak-dis}.
For simplicity, we denote by $\{u_k\}_{k\in\mathbb{N}}$ the
convergent subsequence $\{u_{k_{m_j}}\}_{m\in\mathbb{N}}$, and by
$\{\mathcal{T}_k\}_{k\in\mathbb{N}}$ the subsequence
$\{\mathcal{T}_{k_{m_j}}\}_{m\in\mathbb{N}}$.

We first prove that  limiting eigenpair
$(\lambda_{\infty},u_{\infty})$ is also an eigenpair of
\eqref{weak}. Note that
\begin{eqnarray*}
& &\lambda_{\infty}(u_{\infty},v)-\alpha(\nabla u_{\infty},\nabla v)
- (Vu_{\infty} +\mathcal{N}(u_{\infty}^2)u_{\infty},
v)-\langle\mathbf{R}_k(u_k),v\rangle\\
&=&(\lambda_{\infty}u_{\infty}-\lambda_k
u_k,v)-\alpha(\nabla(u_{\infty}-u_k),\nabla v)\\
& &- (V(u_{\infty}-u_k),v)-(\mathcal{N}(u_{\infty}^2)u_{\infty} -
\mathcal{N}(u^2_{k})u_{k},v) \quad \forall~ v\in H^1_0(\Omega),
\end{eqnarray*}
we obtain from (\ref{n-estimate}) that
\begin{eqnarray}\label{convergence-equ}
& & |\lambda_{\infty}(u_{\infty},v)-\alpha(\nabla u_{\infty},\nabla
v) - (Vu_{\infty} + \mathcal{N}(u_{\infty}^2)u_{\infty},
v)-\langle\mathbf{R}_k(u_k),v\rangle|  \nonumber\\
&\lesssim& \|\nabla(u_{\infty}-u_k)\|_{0,\Omega}\|\nabla
v\|_{0,\Omega}+\|V\|_{0,\Omega}\|u_{\infty}-u_k\|_{0,3,\Omega}
\|v\|_{0,6,\Omega}\nonumber\\
&&+\|u_{\infty}-u_k\|_{1,\Omega}\|v\|_{1,\Omega} +
(\|u_{\infty}-u_k\|_{0,\Omega} +
|\lambda_k-\lambda_{\infty}|)\|v\|_{0,\Omega} \quad\forall ~v\in
H^1_0(\Omega).
\end{eqnarray}
Since $\lambda_k\rightarrow\lambda_{\infty}$ and $u_k\rightarrow
u_{\infty}$ in $H_0^1(\Omega)$, the right hand side of
(\ref{convergence-equ}) tends to zero when $k$ tends to infinity.
Using Lemma \ref{R-limit} and identity
\begin{eqnarray*}
& &\lambda_{\infty}(u_{\infty},v)-\alpha(\nabla u_{\infty},\nabla v)
- (Vu_{\infty} +\mathcal{N}(u_{\infty}^2)u_{\infty},
v)\\
&=&\lambda_{\infty}(u_{\infty},v)-\alpha(\nabla u_{\infty},\nabla v)
- (Vu_{\infty} +\mathcal{N}(u_{\infty}^2)u_{\infty},
v)-\langle\mathbf{R}_k(u_k),v\rangle+
\langle\mathbf{R}_k(u_k),v\rangle,
\end{eqnarray*}
we arrive at
\begin{eqnarray*}
\alpha(\nabla u_{\infty},\nabla v) + (Vu_{\infty} +
\mathcal{N}(u_{\infty}^2)u_{\infty},
v)=\lambda_{\infty}(u_{\infty},v) \quad \forall~ v\in H^1_0(\Omega).
\end{eqnarray*}

Now we  prove that for a sufficiently fine initial mesh, the
limiting eigenfunction $u_{\infty}$ is a ground state solution. Set
\begin{eqnarray*}
\mathcal{W}=\{w\in H_0^1(\Omega)~:~w \mbox{ is an eigenfunction of
\eqref{weak}} \}.
\end{eqnarray*}
Note that $\mathcal{U}\subsetneq\mathcal{W}$, the ground state
solutions in $\mathcal{U}$ minimize energy functional
\eqref{energy}, which is continuous over $H_0^1(\Omega)$, we can
choose a mesh $\mathcal{T}_0$ such that
\begin{eqnarray*}
E_0\equiv E(v) < \min_{w\in \mathcal{W}\setminus\mathcal{U}}E(w)
~~\forall~ v\in \mathcal{U}_0,
\end{eqnarray*}
where the fact \begin{eqnarray*} \lim_{h\rightarrow 0}\inf_{v\in
S_0^h(\Omega)} \|v-w\|_{1,\Omega}=0 \quad \forall~w\in H^1_0(\Omega)
\end{eqnarray*}
is used. Due to $\mathcal {T}_k\subset \mathcal {T}_0$, we have
$E_k\le E_0$ and obtain  $u_{\infty}\in \mathcal{U}$. This completes
the proof.
\end{proof}

If we make a further assumption that
\begin{eqnarray}\label{assumption-convex}
\mathcal{E}''(t)>0\quad \mbox{for}~t\in [0,\infty),
\end{eqnarray}
then energy functional
\begin{eqnarray*}
E(\sqrt{\rho})= \int_{\Omega} \left(\alpha|\nabla \sqrt{\rho}|^2 +
V(x)\rho(x) + \mathcal{E}(\rho(x))\right)dx + \frac{1}{2q}
D_K(\rho^{q},\rho^{q})
\end{eqnarray*}
is strictly convex on  convex set $\{\rho\geq 0:
~\sqrt{\rho}\in\mathcal{A}\}$ and hence there exists a unique
minimizer of \eqref{energy} in  admissible class $\mathcal{A}$. Note
that the minimizer  of \eqref{energy} in $V_k$ is unique when
initial mesh $\mathcal{T}_0$ is fine enough (c.f., e.g.,
\cite{zhou07}), we have

\begin{corollary}
Assume that the hypothesis of Theorem \ref{theorem-convergence} and
\eqref{assumption-convex} are satisfied. If $(\lambda,u)\in
\mathbb{R} \times \mathcal {A}$ is the ground state solution of
\eqref{weak} and $(\lambda_k,u_k)\in \mathbb{R} \times V_k$ is the
ground state solution of \eqref{weak-dis}, then
\begin{eqnarray*}
\lim_{k\rightarrow \infty}\|u_k-u\|=0,\\[1ex]
\lim_{k\rightarrow\infty}\|\lambda_k-\lambda\|=0.
\end{eqnarray*}
\end{corollary}

\section{Numerical examples}\label{sec-example}\setcounter{equation}{0}
In this section, we will report on some numerical experiments for
both linear finite elements and quadratic finite elements in three
dimensions to illustrate the convergence of adaptive finite element
approximations.

Our numerical computations are carried out on LSSC-II in the State
Key Laboratory of Scientific and Engineering Computing, Chinese
Academy of Sciences, and our codes are based on the toolbox PHG of
the laboratory. All of the computational results are given in atomic
unit (a.u.).

{\bf Example 1.} Consider the ground state solution of GPE for BEC
with a harmonic oscillator potential
\begin{eqnarray*}
V(x,y,z)=\frac{1}{2}(\gamma_x^2x^2+\gamma_y^2y^2+\gamma_z^2z^2),
\end{eqnarray*}
where $\gamma_x=1,\gamma_y=2,\gamma_z=4$. We solve the following
nonlinear problem: Find $(\lambda,u)\in \mathbb{R}\times
H^1_0(\Omega)$ such that $\|u\|_{0,\Omega}=1$ and
\begin{eqnarray*}
\left\{\ \begin{array}{rclc} \left(-\frac{1}{2}\Delta + V +
\beta |u|^2\right) u &=& \lambda u & \mbox{in}~\Omega,\\[1ex]
u&=&0&\mbox{on}~\partial\Omega,
\end{array} \right.
\end{eqnarray*}
where $\beta=200$ and $\Omega=[-8,8]\times [-6,6]\times [-4,4]$.

\begin{figure}[ht]
\centering
\includegraphics[width=6.3cm]{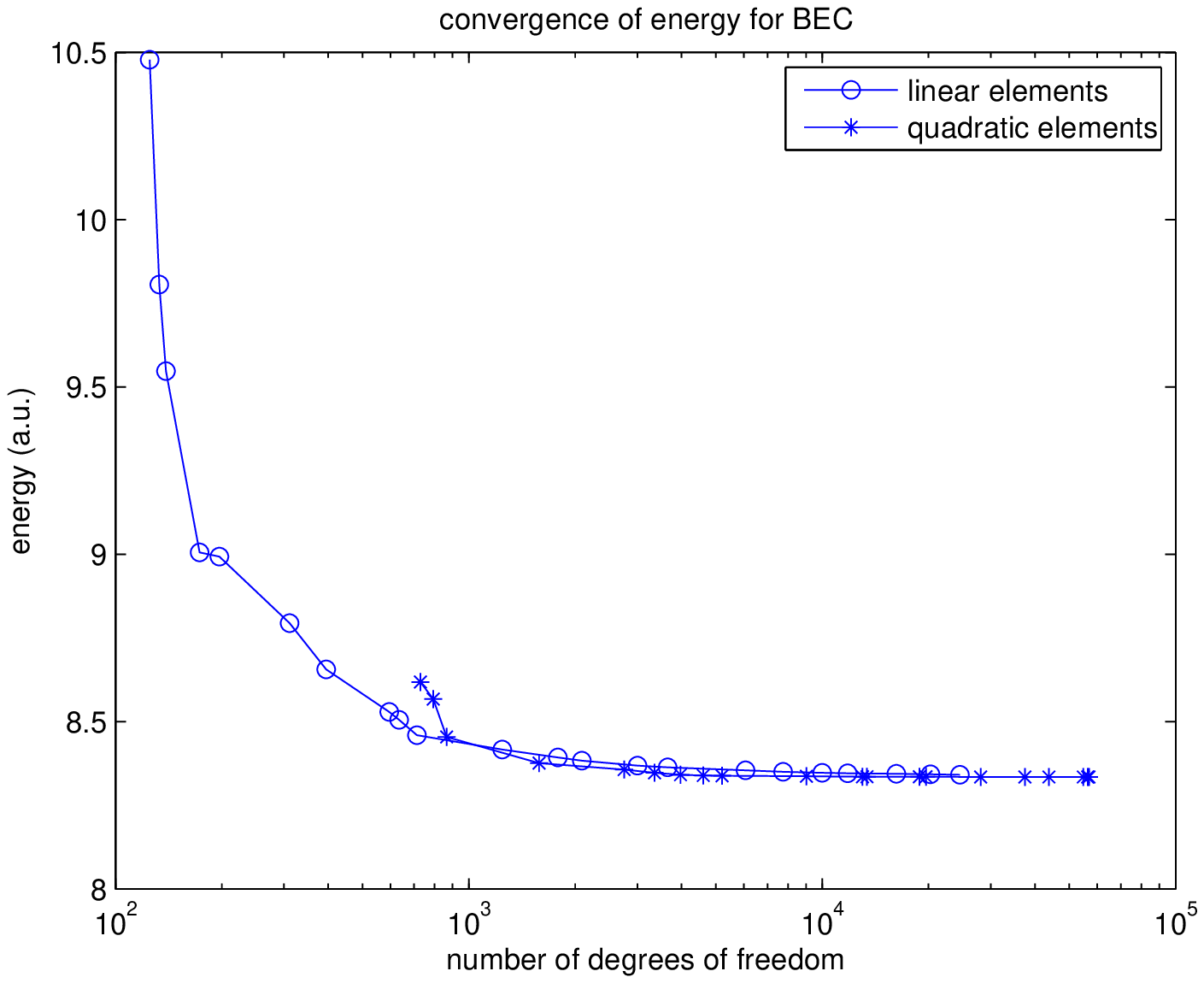}
\includegraphics[width=6.3cm]{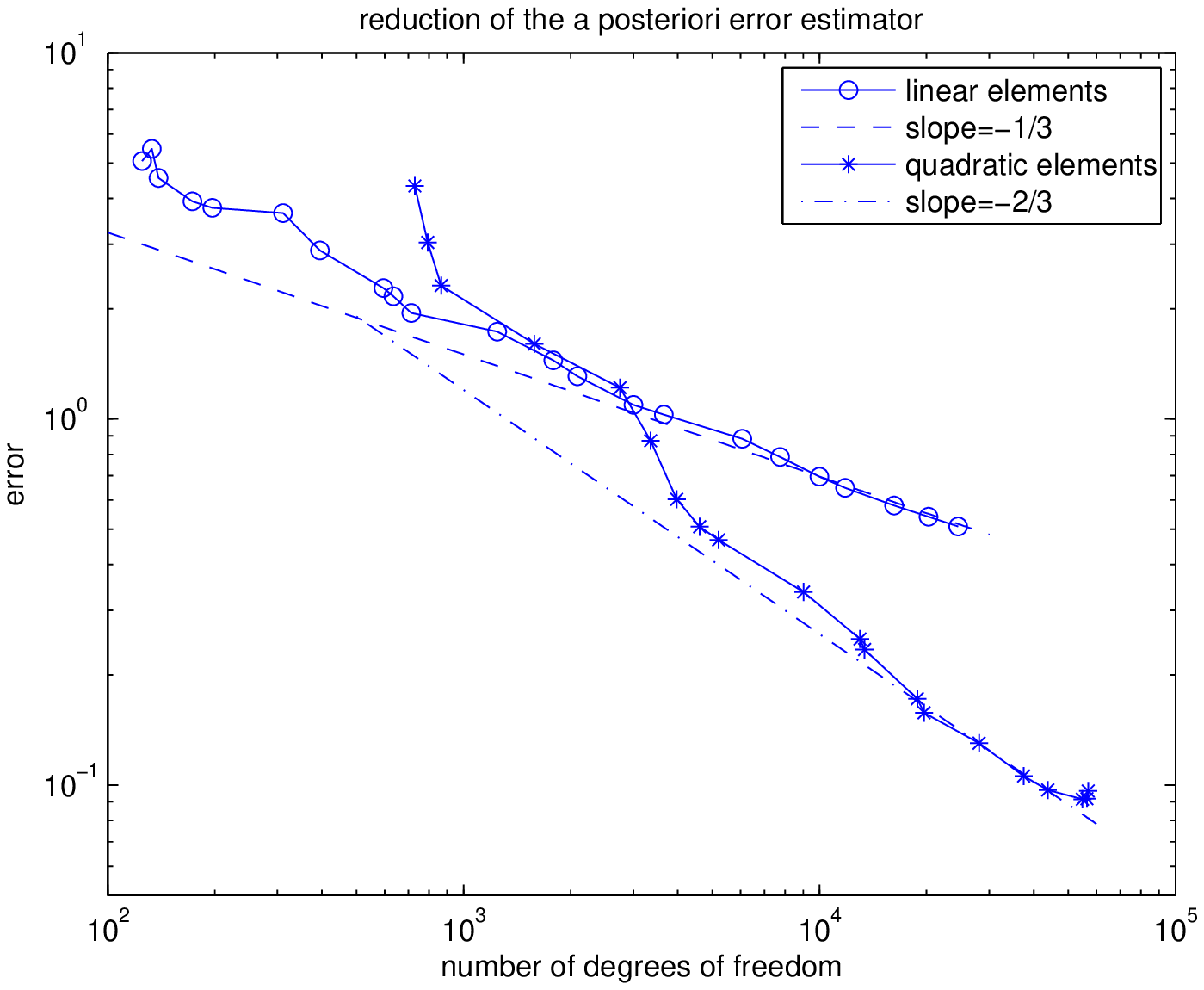}
\caption{Left: Convergence curves of energy for BEC. Right:
Reduction of the a posteriori error estimators using  linear and
quadratic elements.}\label{example1-1}
\end{figure}

\begin{figure}[ht]
\centering
\includegraphics[width=5.5cm]{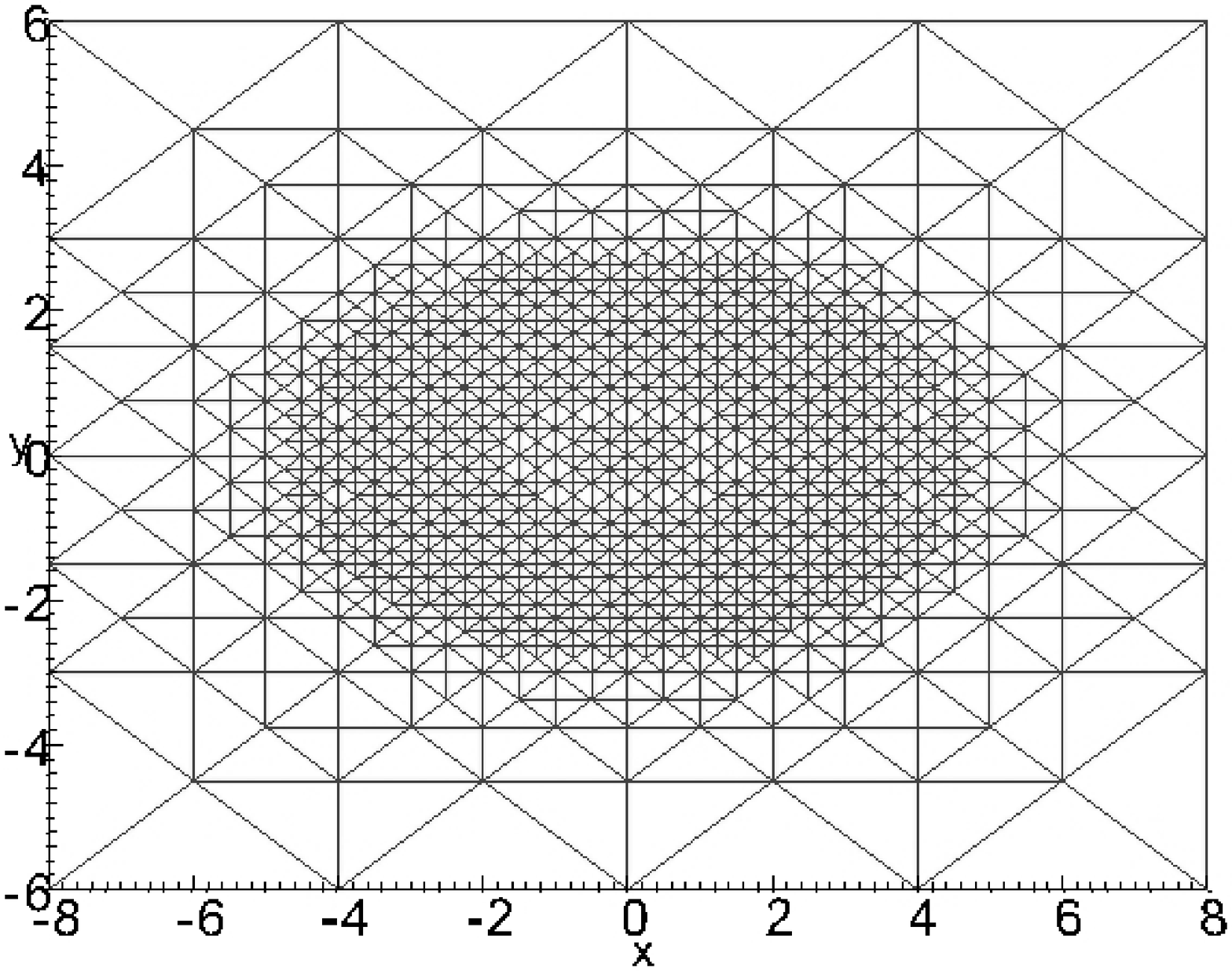}\hskip 1.0cm
\includegraphics[width=5.5cm]{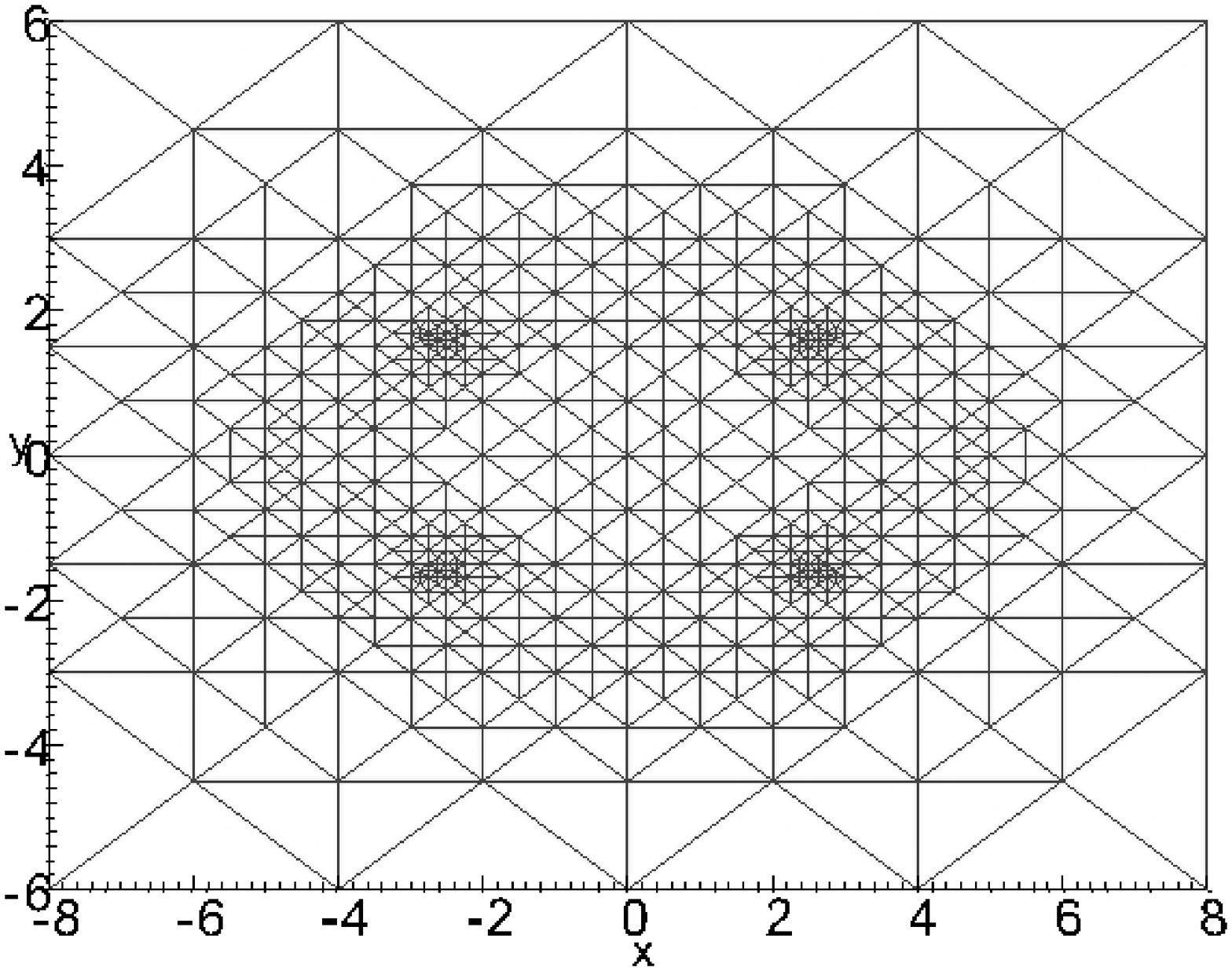}
\caption{The cross-sections on $z=0$ of adaptive meshes using linear
(left) and quadratic (right) elements.}\label{example1-2}
\end{figure}

The convergence of  energies and the reduction of the a posteriori
error estimators are presented in Figure \ref{example1-1}, which
support our theory and that the a posteriori error estimators are
efficient. Some cross-sections of the adaptively refined meshes
constructed by the a posteriori error indicators are displayed in
Figure \ref{example1-2}.


In the next two examples, we shall carry out the ground state energy
calculations of atomic and molecular systems based on TFW type
orbital-free models. The nonlinear term is given by
\begin{eqnarray*}\label{f--00}
\mathcal{N}(u^2) = \int_{\mathbb{R}^3}\frac{u^2(y)}{|\cdot-y|}dy +
\frac{5}{3} C_{TF}u^{4/3}+ \upsilon_{xc}(u^2),
\end{eqnarray*}
where $C_{TF}=\frac{3}{10}(3\pi^2)^{\frac{2}{3}}$ and $v_{xc}$ is
the exchange-correction potential. The exchange-correction potential
used in our computation is chosen as
\begin{eqnarray}\label{f-2}
\upsilon_{xc}(\rho)=\upsilon_{x}^{LDA}(\rho)+\upsilon_{c}^{LDA}(\rho),
\end{eqnarray}
where
\begin{eqnarray*}
\upsilon_{x}^{LDA}(\rho)=-\left(\frac{3}{\pi}\right)^{1/3}\rho^{1/3},
\end{eqnarray*}
\begin{eqnarray*}
\upsilon_{c}^{LDA}(\rho)=\left\{
\begin{array}{cc}
0.0311\ln r_s-0.0584+0.0013r_s\ln r_s-0.0084r_s &\mbox{if}~r_s<1\\[0.2cm]
-(0.1423+0.0633r_s+0.1748\sqrt{r_s})/(1+1.0529\sqrt{r_s}+0.3334r_s)^2
&\mbox{if}~r_s\geq 1
\end{array}\right.
\end{eqnarray*}
and $\displaystyle r_s=\left(\frac{3}{4\pi\rho}\right)^{1/3}$.

{\bf Example 2.} Consider the TFW type orbital-free model for helium
atoms. The external electrostatic potential is  $V(x)=-2/|x|.$ Then
we have the following nonlinear problem: Find $(\lambda,u)\in
\mathbb{R}\times H^1_0(\Omega)$ such that $\|u\|^2_{0,\Omega}=2$ and
\begin{eqnarray*}
\left\{\ \begin{array}{rclc}\displaystyle -\frac{1}{10}\Delta
u-\frac{2}{|x|}u+u\int_{\Omega}\frac{|u(y)|^2}{|x-y|}dy
+\frac{5}{3}C_{TF}u^{7/3}+v_{xc}(u^2)u &=& \lambda u
& \mbox{in}~\Omega,\\[1ex] u&=&0&\mbox{on}~\partial\Omega,
\end{array} \right.
\end{eqnarray*}
where $\Omega =(-5.0,5.0)^3$.

\begin{figure}[ht]
\centering
\includegraphics[width=6.3cm]{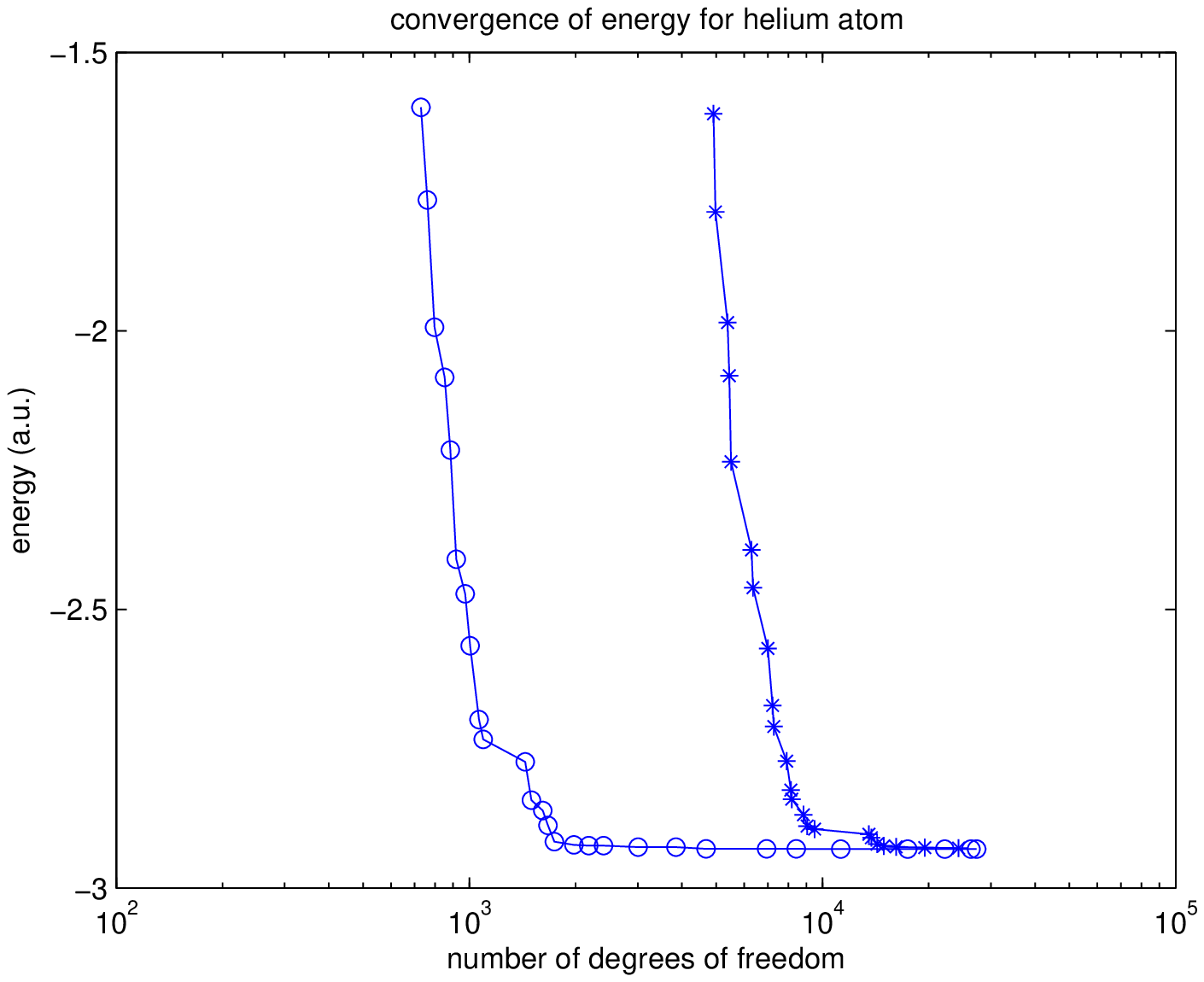}
\includegraphics[width=6.3cm]{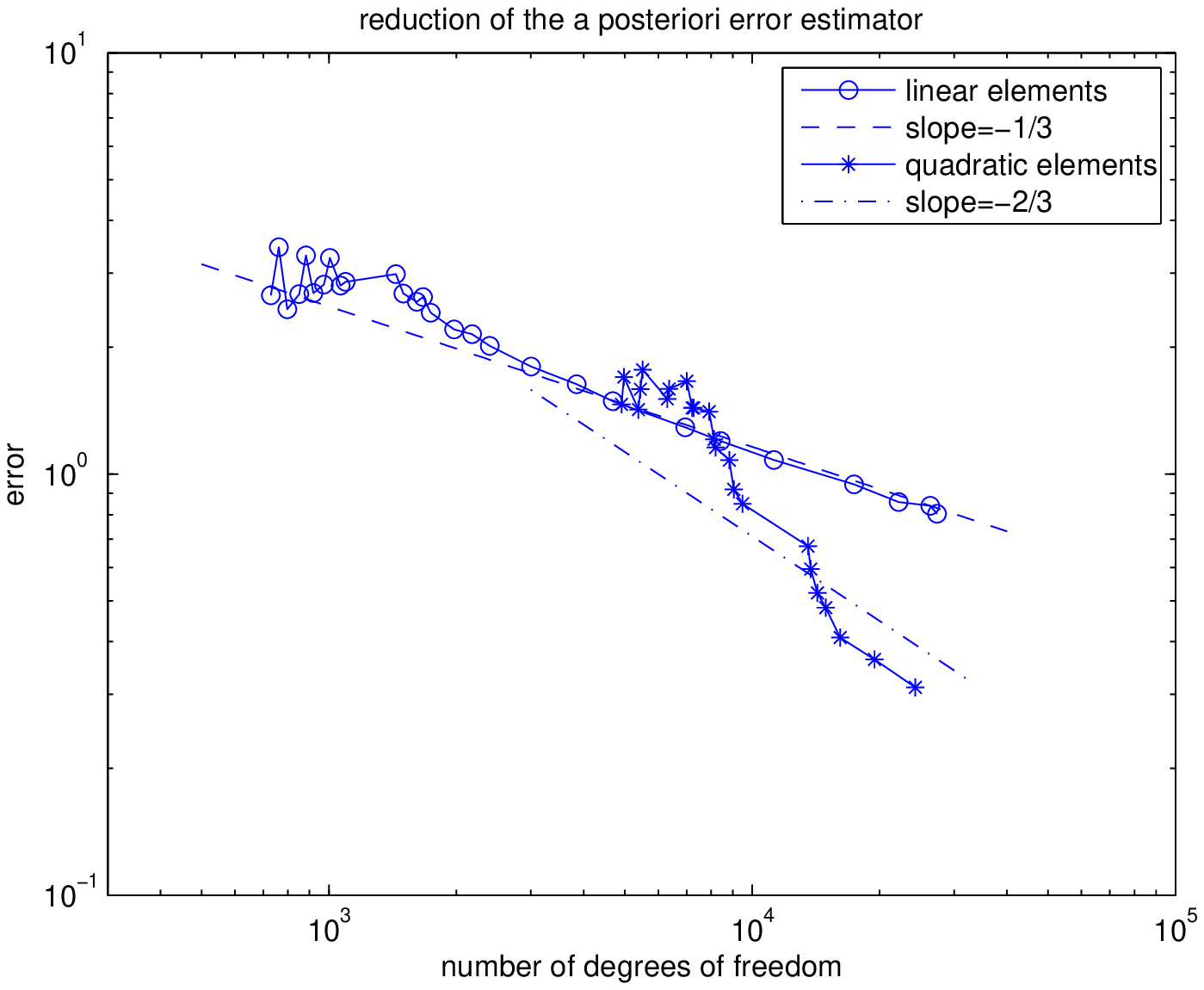}
\caption{Left: Convergence curves of energy for the helium atom.
Right: Reduction of the a posteriori error estimators using linear
and quadratic elements.}\label{example2-1}
\end{figure}

\begin{figure}[ht]
\centering
\includegraphics[width=6.3cm]{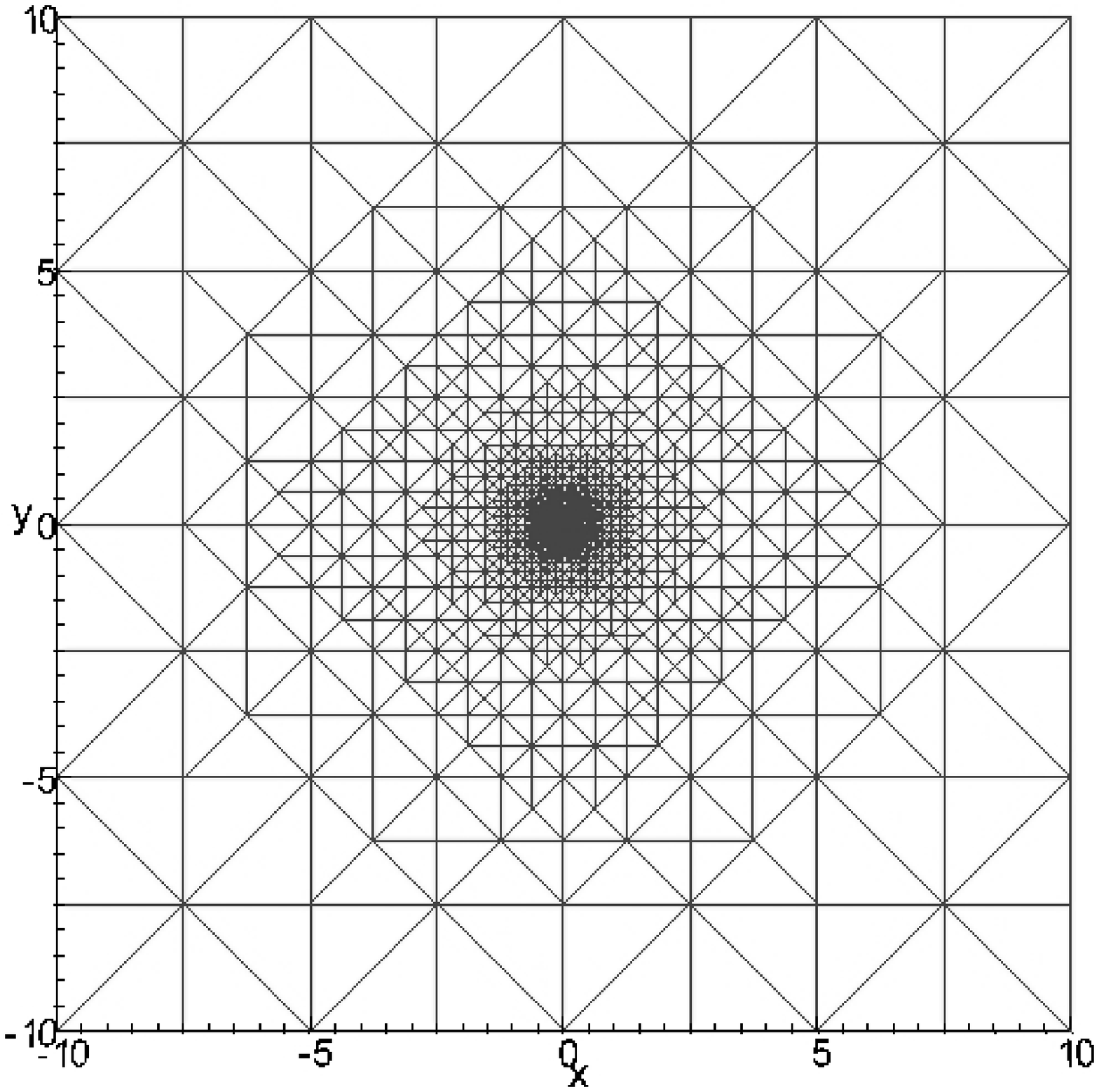}\hskip 0.6cm
\includegraphics[width=6.3cm]{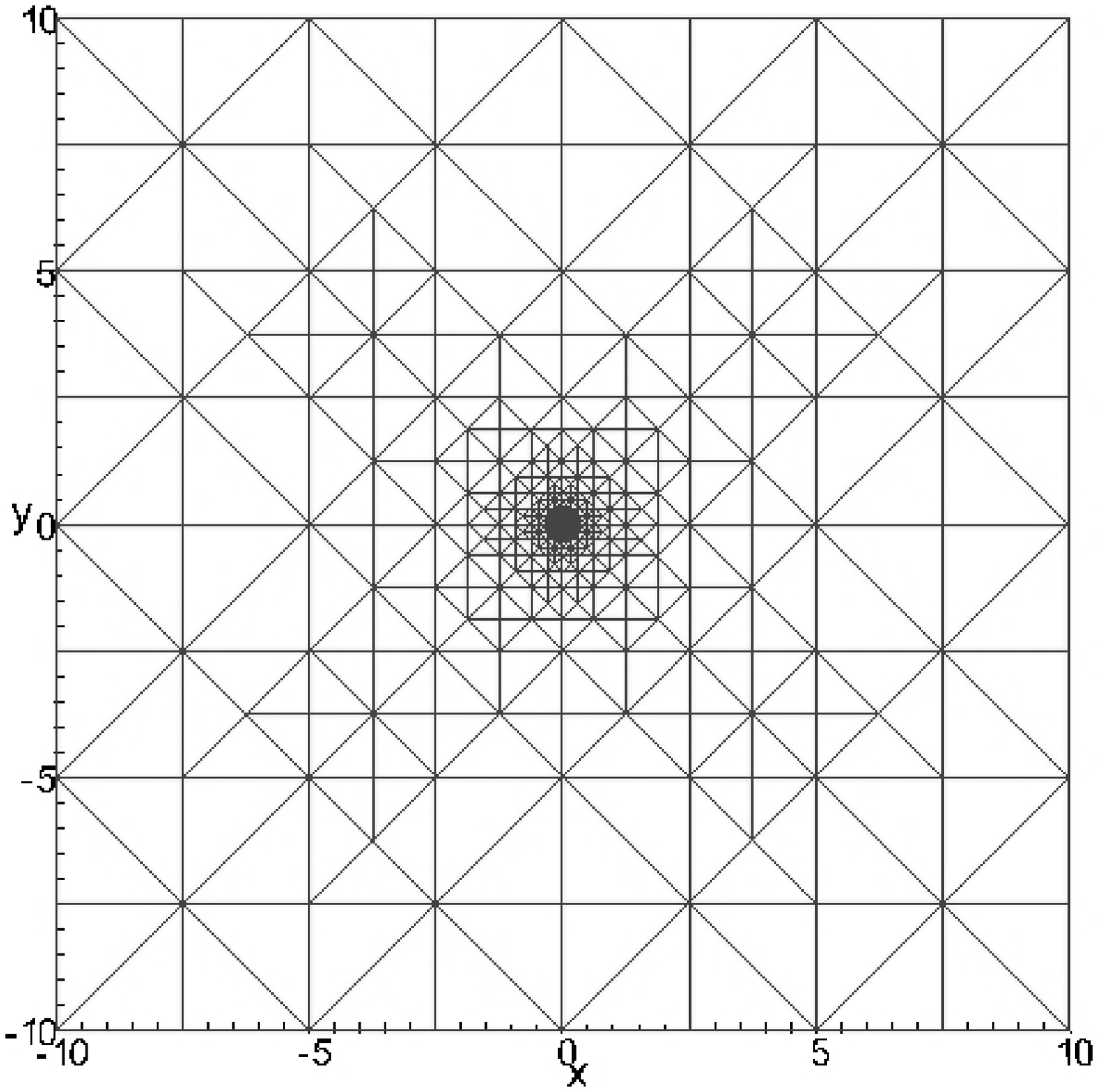}
\caption{The cross-sections on $z=0$ of adaptive meshes using linear
(left) and quadratic (right) elements.}\label{example2-2}
\end{figure}

The convergence of energies and the reduction of the a posteriori
error estimators are shown in Figure \ref{example2-1}, which support
our theory. The cross-sections of the adaptive meshes are displayed
in Figure \ref{example2-2}, from which we observe that more refined
meshes (nodes) appear in the area where the nuclei are located.


{\bf Example 3.} Finally, we consider an aluminum cluster in the
face centered cubic lattice consisting of $3\times 3\times 3$ unit
cells with 172 aluminium atoms, where the GHN pseudopotential
\cite{goodwin-etal90} is used. We solve the following nonlinear
problem: Find $(\lambda,u)\in \mathbb{R}\times H^1_0(\Omega)$ such
that $\|u\|^2_{0,\Omega}=172$ and
\begin{eqnarray*}
\left\{\ \begin{array}{rclc}\displaystyle -\frac{1}{10}\Delta u +
V_{pseu}^{GHN} u + u\int_{\Omega}\frac{|u(y)|^2}{|x-y|}dy +
\frac{5}{3}C_{TF}u^{7/3} + v_{xc}(u^2)u &=& \lambda u &
\mbox{in}~\Omega,\\[1ex] u &=& 0&\mbox{on}~\partial\Omega,
\end{array} \right.
\end{eqnarray*}
where $\Omega=(-25.0,25.0)^3$.

\begin{figure}[ht]
\centering
\includegraphics[width=6.3cm]{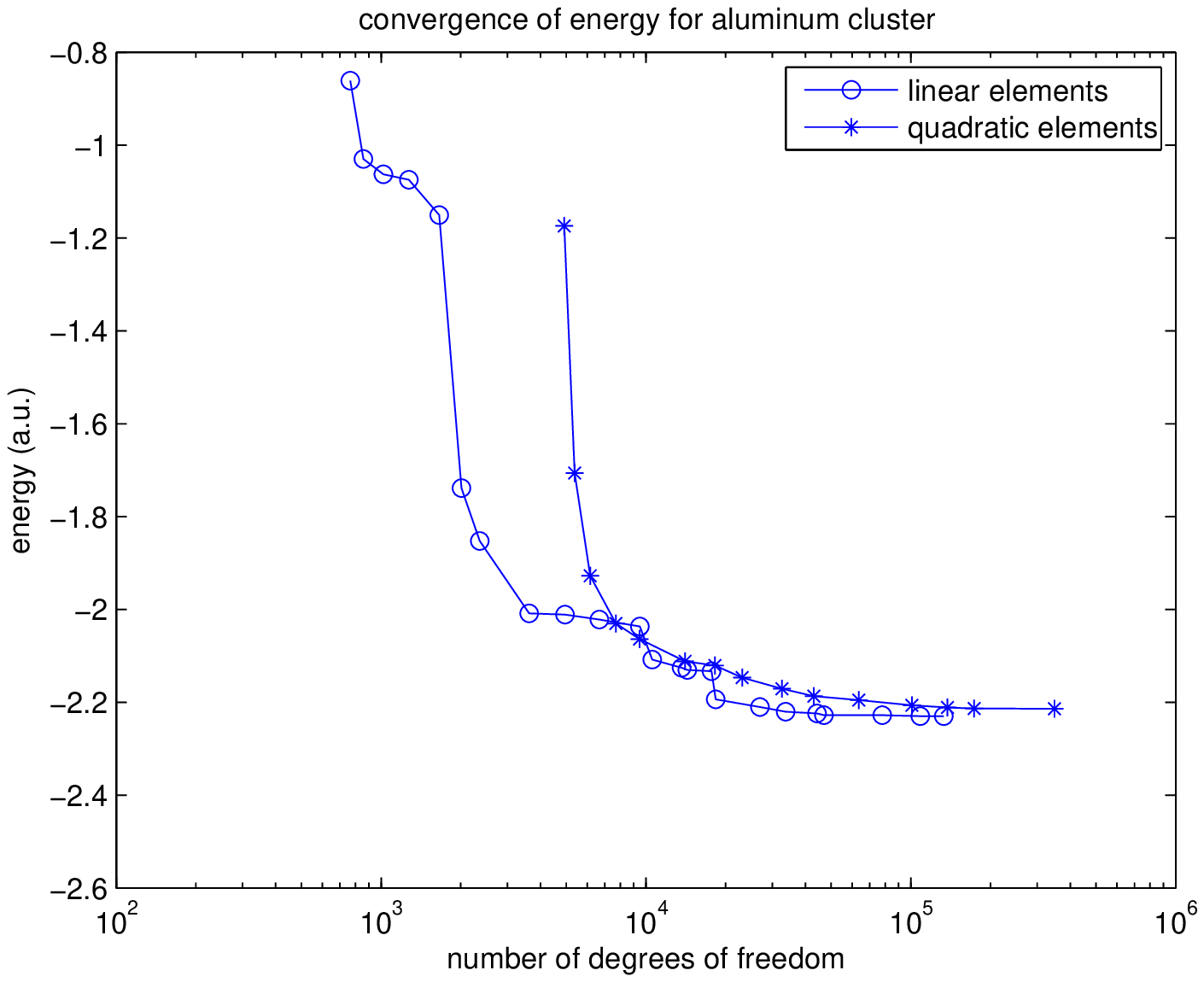}
\includegraphics[width=6.3cm]{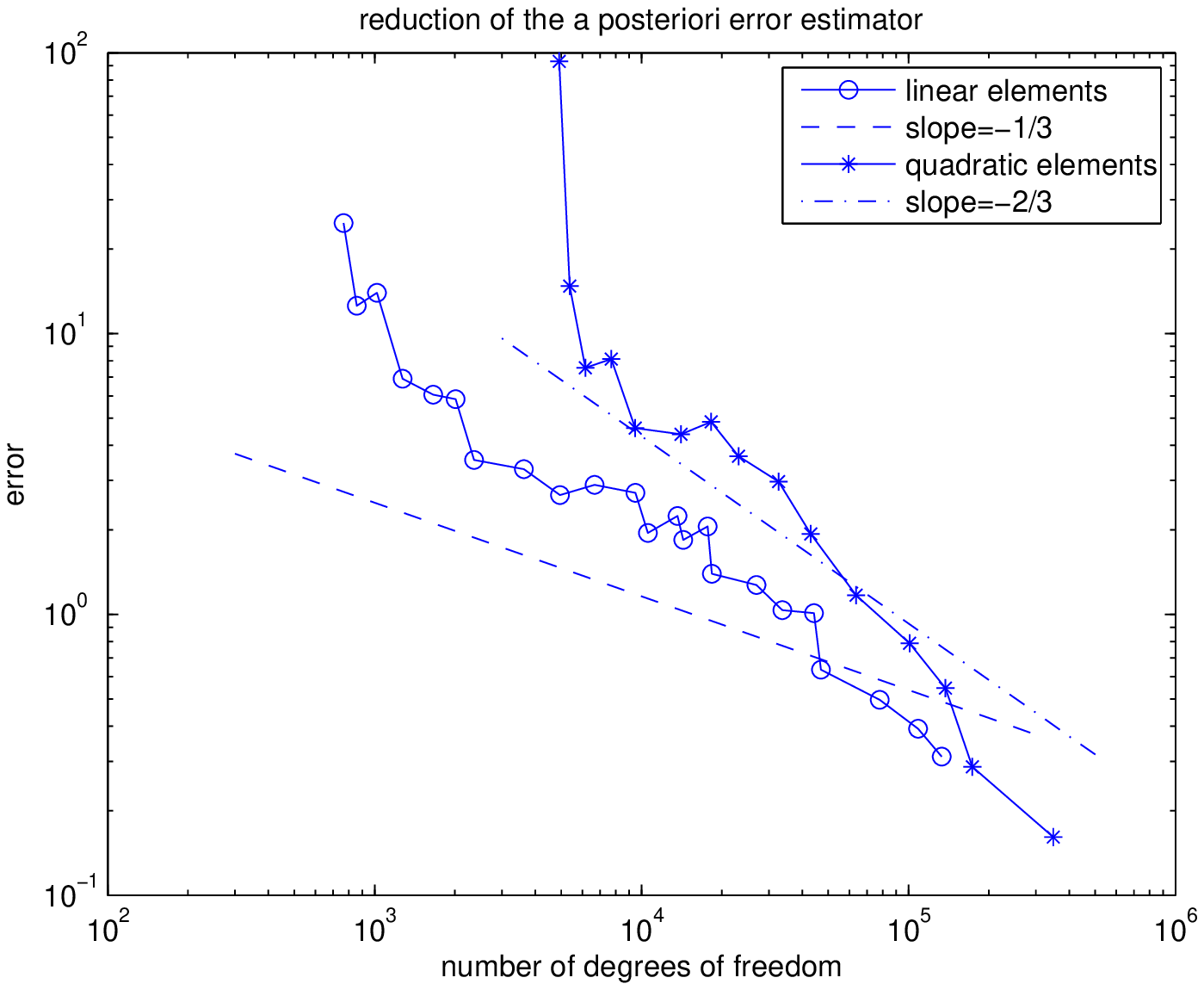}
\caption{Left: Convergence curves of energy for the aluminium
cluster in FCC lattice. Right: Reduction of the a posteriori error
estimators using linear and quadratic elements.}\label{example3-1}
\end{figure}

The convergence of  energies and the reduction of a posteriori error
estimators are shown in Figure \ref{example3-1}. The cross-sections
of the adaptive meshes are displayed in Figure \ref{example2-2}. We
observe that with the a posteriori error estimators, the refinement
is carried out automatically at the regions where the computed
functions vary rapidly, especially near the nuclei. As a result, the
computational accuracy can be controlled efficiently and the
computational cost is reduced significantly.

\begin{figure}[ht]
\centering
\includegraphics[width=6.3cm]{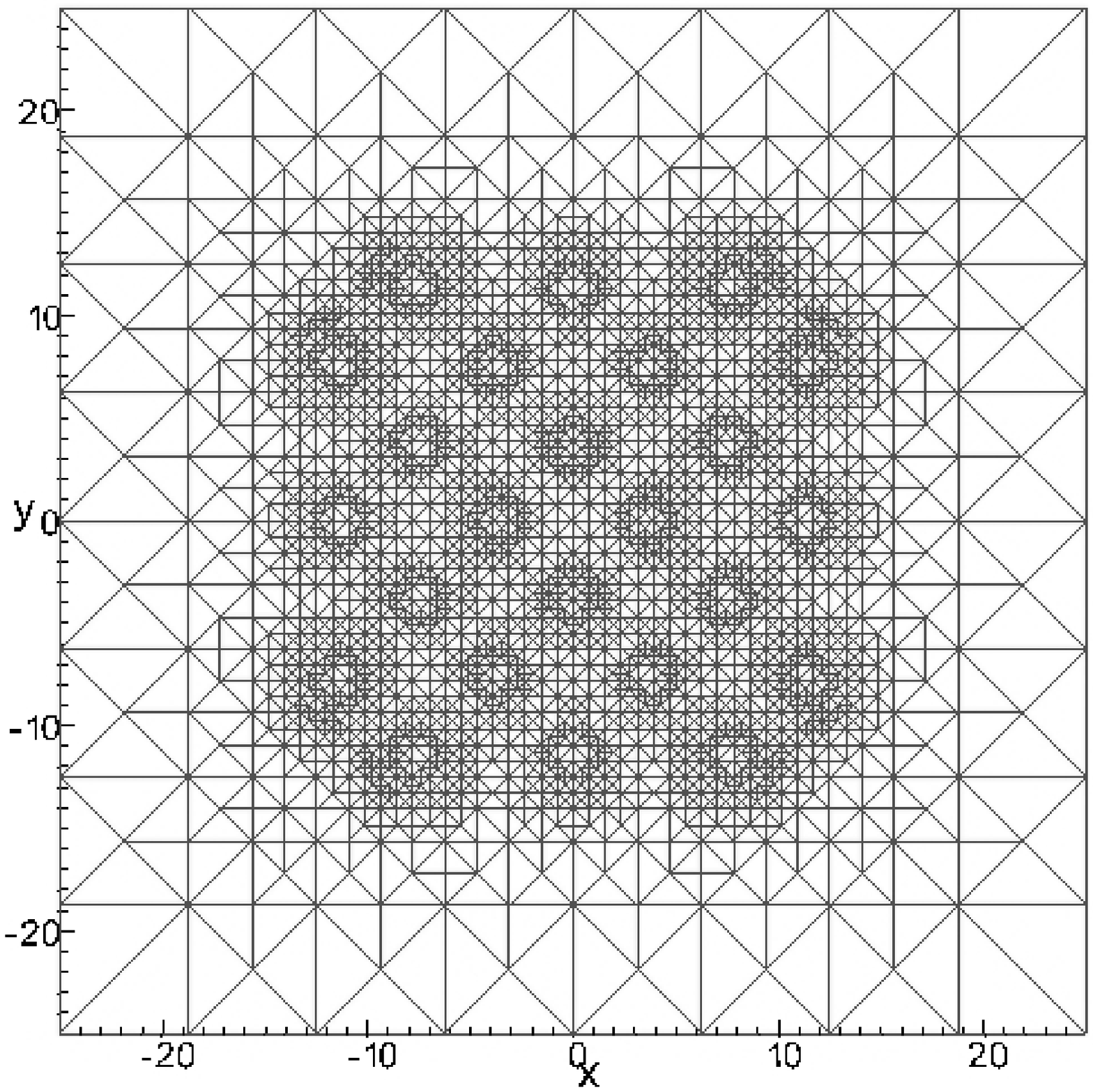}\hskip 0.6cm
\includegraphics[width=6.3cm]{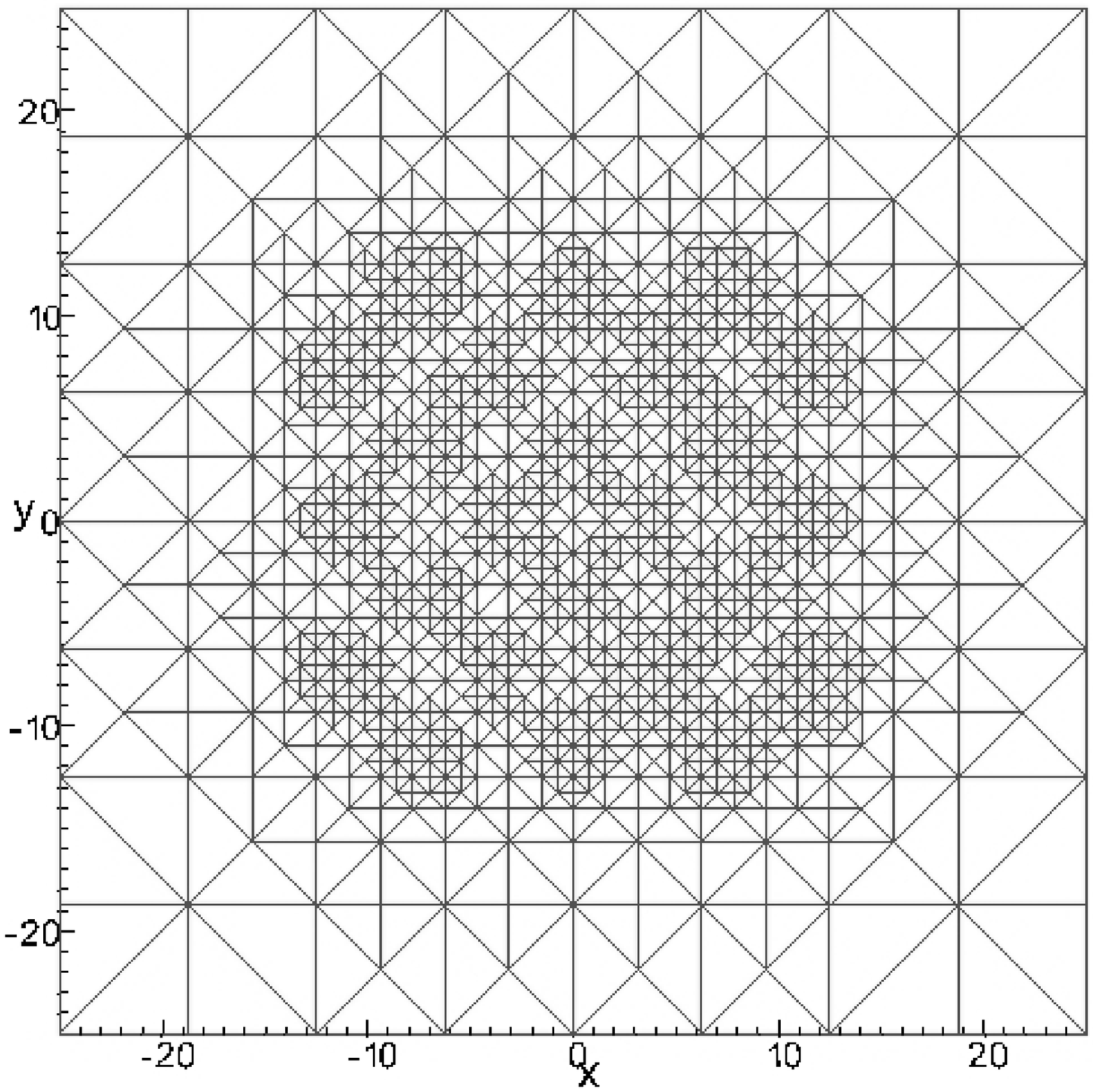}
\caption{The cross-sections on $z=0$ of adaptive meshes using linear
and quadratic finite elements.}\label{example3-2}
\end{figure}

\section{Concluding remarks}\setcounter{equation}{0}
We have analyzed adaptive finite element approximations for ground
state solutions of a class of nonlinear eigenvalue problems. We have
proved that the adaptive finite element loop produces a sequence of
approximations that converge to the set of exact ground state
solutions. This result covers many mathematical models of practical
interest, for instance, the Bose-Einstein condensation, the TFW
model in the orbital-free density functional theory, and
Schr\"{o}dinger-Newton equations in the quantum state reduction
\cite{chen-gong-zhou09,harrison03,penrose96} where the integration
kernel $K$ is negative. We have also applied  adaptive finite
element discretizations to micro-structure of matter calculations,
which support our theory. It is shown by Figure \ref{example1-1},
Figure \ref{example2-1}, and Figure \ref{example3-1} that we may
have some convergence rates of adaptive finite element
approximations. Indeed, it is our on-going work to study the optimal
complexity of adaptive finite element approximations for such
nonlinear eigenvalue problems, which requires a new technical tool
and will be addressed elsewhere \cite{chen-he-zhou10}.\vskip 0.2cm

{\sc Acknowledgements.} The authors would like to thank Dr. Xiaoying
Dai,   Prof. Lihua Shen, and Dr. Dier Zhang for their stimulating
discussions and fruitful cooperations on electronic structure
computations that have motivated this work. The authors are grateful
to Prof. Linbo Zhang and Mr. Hui Liu for their assistance on
numerical computations.

\section*{Appendix}
\renewcommand{\theequation}{A.\arabic{equation}}
\renewcommand{\thetheorem}{\arabic{section}.\arabic{theorem}}
\renewcommand{\thelemma}{A.\arabic{lemma}}
\renewcommand{\theproposition}{A.\arabic{proposition}}
\renewcommand{\thesection}{A}
\setcounter{equation}{0} We may follow the framework in
\cite{verfurth96} to derive an upper bound and a lower bound of the
a posteriori error estimate. Let
\begin{gather*}
X~=~Y~=~H_0^1(\Omega),\\
\|\cdot\|_X~=~\|\cdot\|_Y~=~\|\cdot\|_{1,\Omega},\\
\langle F(\lambda,u),(\mu,v) \rangle = \int_{\Omega} \big( \alpha
\nabla u\nabla v+Vuv+\mathcal{N}(u^2)uv-\lambda uv \big)
+\mu\big(\int_{\Omega}u^2-Z \big).
\end{gather*}
It is seen that
 (\ref{weak}) may be written
as
\begin{eqnarray}\label{eq-abstract}
F(\lambda,u)=0
\end{eqnarray}
and $F\in C^1(\mathbb{R}\times X,\mathbb{R}\times Y^*)$, where $Y^*$
is the dual space of $Y$. A solution $(\lambda,u)$ of
\eqref{eq-abstract} is said to be regular  if  the following
equation
$$
DF(\lambda,u)\cdot(\mu,v)=(\kappa,f)
$$
is uniquely solvable in $\mathbb{R}\times X$ for each $(\kappa,f)\in
\mathbb{R}\times Y^*$, where $DF(\lambda,u):~\mathbb{R}\times
X\rightarrow\mathbb{R}\times Y^*$ is the Fr\'{e}chet derivative of
$F$ at $(\lambda,u)$. It is seen that $(\lambda,u)\in
\mathbb{R}\times H^1_0(\Omega)$ is a regular solution of
(\ref{eq-abstract}) if
\begin{eqnarray}\label{eq-regular}
\langle (E''(u)-\lambda)v,v \rangle_{Y^*,X} \geq \gamma\|\nabla
v\|^2_{0,\Omega}\quad \forall~ v\in H_0^1(\Omega)\cap u^{\bot}
\end{eqnarray}
is true for some constant $\gamma >0$ (see, e.g.,
\cite{langwallner-ortner-suli09}), where
\begin{eqnarray*}
\langle (E''(u)-\lambda)v,w \rangle &=& \alpha(\nabla v,\nabla w) +
((V+\mathcal{N}(u^2)-\lambda)v, w) +
2(\mathcal{N}_1'(u^2)u^2v,w) \\[1ex] && +
2qD_K(u^{2q-1}v,u^{2q-1}w)+2(q-1)(\mathcal{N}_2(u^2)v,w)
\end{eqnarray*}
and $$ u^{\bot}=\{w\in L^2(\Omega): (u,w)=0\}.$$ It has been proved
that (\ref{eq-regular}) is satisfied by some special TFW models that
are of convex functional (see, e.g.,
\cite{cances-chakir-maday09a,cances-chakir-maday09b}).

Let
\begin{gather*}
X_h~=~Y_h~=~S_0^{h}(\Omega)~\end{gather*} and define
\begin{gather*} \langle
F_h(\lambda_h,u_h),(\mu,v) \rangle ~=~ \langle
F(\lambda_h,u_h),(\mu,v) \rangle \quad
\forall~(\mu,v)\in\mathbb{R}\times S_0^h(\Omega).
\end{gather*}
We see that  $F_h\in C(\mathbb{R}\times X_h,\mathbb{R}\times Y_h^*)$
and is an approximation of $F$.  Obviously, finite element
eigenvalue problem (\ref{weak-dis}) is equivalent to
\begin{eqnarray*}\label{eq-abstract-dis}
F_h(\lambda_h,u_h)=0
\end{eqnarray*}
 and $(\lambda_h,u_h)$ is an approximate
solution of (\ref{eq-abstract}).

The following proposition in \cite[Section 2.1]{verfurth96} yields a
posteriori error estimates  in the neighborhood of $(\lambda,u)$
that satisfies (\ref{eq-abstract}).
\begin{proposition}
Let $(\lambda,u)$ be a regular solution of (\ref{eq-abstract}). If
$DF$ is the derivative of $F$ and $DF$ is Lipschitz continuous at
$(\lambda,u)$, then the following estimate holds for all
$(\lambda_h,u_h)$ sufficiently close to this solution:
\begin{eqnarray}\label{estimate-2.1}
\|F(\lambda_h,u_h)\|_{\mathbb{R} \times Y^*} \lesssim
|\lambda_h-\lambda|+\|u_h-u\|_X \lesssim
\|F(\lambda_h,u_h)\|_{\mathbb{R} \times Y^*}.
\end{eqnarray}
\end{proposition}

It is shown from (\ref{estimate-2.1}) that
$\|F(\lambda_h,u_h)\|_{\mathbb{R} \times Y^*}$  is a posteriori
error estimator. When we  apply the general approach in
\cite[Sections 3.3-3.4]{verfurth96} to
\begin{gather*}
\underline{a}(x,u,\nabla u)=\alpha\nabla u \quad\mbox{and}\quad
b(x,u,\nabla u)=\lambda u-Vu-\mathcal{N}(u^2)u
\end{gather*}
by taking
\begin{gather*}
\underline{a}_h(x,u_h,\nabla u_h)=\alpha\nabla u_h, \quad
b_h(x,u_h,\nabla u_h)=\overline{\lambda u_h} - \overline{Vu_h} -
\overline{\mathcal{N}(u_h^2)u_h},
\end{gather*}
and
\begin{gather*}
\eta_T = \eta_h(u_h,T),\quad \varepsilon_T = {\rm osc}_h(u_h,T),
\end{gather*}
we then gives a  proof of Theorem \ref{theorem-bound}.

\end{document}